\documentclass[10pt]{article}

\newcommand{\mylabel}
{\label}

\usepackage{amssymb}
\usepackage{amsthm}
\usepackage{amsmath}
\usepackage{amscd}
\usepackage{multicol}

\allowdisplaybreaks

\newcommand{\R}{\mathbf{R}}
\newcommand{\C}{\mathbf{C}}
\newcommand{\eps}{\varepsilon}

\newcommand{\me}{e}
\newcommand{\mi}{i}

\newcommand{\dif}{d}
\newcommand{\Dif}{\mathrm{D}}

\renewcommand{\det}{\mathrm{det}}

\newcommand{\papprox}{\tilde S_{a}}
\newcommand{\cobg}{C^{0,\alpha}_{\gamma-2}}
\newcommand{\ctbg}{C^{2,\alpha}_\gamma}
\newcommand{\ckbg}{C^{k,\alpha}_\gamma}
\newcommand{\Bkbg}{\mathcal{B}^{k,\alpha,\gamma}_a}
\newcommand{\Btbg}{\mathcal{B}^{2,\alpha, \gamma}_a}
\newcommand{\Bobg}{\mathcal{B}^{0,\alpha,\gamma-2}_a}

\newtheorem{thm}{Theorem}
\newtheorem{lemma}[thm]{Lemma}
\newtheorem{cor}[thm]{Corollary}
\newtheorem{prop}[thm]{Proposition}
\newtheorem*{nonumthm}{Theorem}

\theoremstyle{definition}
\newtheorem{defn}[thm]{Definition}

\long\def\symbolfootnote[#1]#2{\begingroup%
\def\thefootnote{\fnsymbol{footnote}}\footnote[#1]{#2}\endgroup}

\begin{document}

\title{Doubling Constant Mean Curvature Tori in $S^3$}

\author{
\begin{minipage}{3.125in}
    \begin{center}
        Adrian Butscher \\ University of Toronto at Scarborough \\ email: \ttfamily butscher@utsc.utoronto.ca
    \end{center}
\end{minipage}
\begin{minipage}{3.125in}
\begin{center}
        Frank Pacard \\ Universit\'e de Paris XII \\ email: \ttfamily pacard@univ-paris12.fr
    \end{center}
\end{minipage}\\
\rule{1ex}{0ex}
}

\maketitle

\begin{abstract}
The Clifford tori in $S^3$ constitue a one-parameter family of flat, two-dimensional, constant mean curvature (CMC) submanifolds. This paper demonstrates that new, topologically non-trivial CMC surfaces resembling a pair of neighbouring Clifford tori connected at a sub-lattice consisting of at least two points by small catenoidal bridges can be constructed by perturbative PDE methods. That is, one can create a submanifold that has almost everywhere constant mean curvature by gluing a re-scaled catenoid into the neighbourhood of each point of a sub-lattice of the Clifford torus; and then one can show that a constant mean curvature perturbation of this submanifold does exist.
\end{abstract}

\renewcommand{\baselinestretch}{1.25}
\normalsize

\section{Introduction and Statement of Results}

\paragraph*{CMC surfaces.}\symbolfootnote[0]{MSC Classification numbers: 53A10, 58J10}
\ \hspace{-3ex}A constant mean curvature (CMC) surface $\Sigma$ contained in an ambient Riemannian manifold $X$ has the property that its mean curvature with respect to the induced metric is constant.  This property ensures that the surface area of $\Sigma$ is a critical value of the area functional for surfaces of $X$ subject to an  enclosed-volume constraint. Constant mean curvature surfaces have been objects of great interest since the beginnings of modern differential geometry.  Classical examples of non-trivial CMC surfaces in $\R^3$ are the sphere, the cylinder and the Delaunay surfaces, and for a long while these were the only known CMC surfaces.  In fact, a result of Alexandrov \cite{alexandrov} states that the only compact, connected, embedded CMC surfaces in $\R^3$ are spheres.

In recent decades, the theory of CMC surfaces in $\R^3$ has  progressed considerably. In 1986, Wente discovered a family of compact, immersed CMC tori \cite{wente}; these have been thoroughly studied also in \cite{pinkallsterling}.  Since then, several parallel sequences of ideas have led to a profusion of new CMC surfaces.  First, the techniques used by Wente have culminated in a representation for CMC surfaces in $\R^3$ akin to the classically-known Weierstra\ss\ representation of minimal surfaces in which a harmonic but not anti-conformal map from a Riemann surface to the unit sphere becomes the Gau\ss\ map of a CMC immersion into $\R^3$ from which the immersion can be determined \cite{dpw,kenmotsu}.  Amongst the many examples of its use are the periodic CMC surfaces of $\R^3$ constructed by Ritor\'e \cite{ritore}.  Second, Kapouleas pioneered the use of geometric partial differential equations to construct many new CMC surfaces: e.g.~compact genus 2 surfaces by fusing Wente tori \cite{kapouleas5}; and compact surfaces of higher genus and non-compact surfaces with arbitrary numbers of ends by gluing together spheres and Delaunay surfaces \cite{kapouleas3,kapouleas2}. Kapouleas' discoveries have since been complemented by much research into gluing and other analysis-based constructions that can be performed in the class of CMC surfaces, most notably in the work of Mazzeo, Pacard and others \cite{jlelipacard,mazzeopacard1,mazzeopacard2,mazzeopacardpollack,mazzeo}. Finally, an idea of Lawson \cite{lawson} coupled with a Schwartz reflection principle allows CMC surfaces of $\R^3$ to be constructed via an associated minimal surface of the 3-sphere.  A great many CMC surfaces containing periodic regions have been constructed by Karcher and Gro{\ss}e-Brauckmann in this way \cite{grossebrauckmann,karcher}.

The corresponding picture amongst CMC surfaces of the sphere $S^3$ is not as rich.  The
classically known examples are the spheres obtained from intersecting $S^3$ with hyperplanes, and the so-called Clifford tori $C_a$, given by
\[ 
C_a : =  \left\{(z_1, z_2) \in \C^2 \, : \, |z_1| =\sqrt{\frac{_{1+a}}{^2}} \quad |z_2| = \sqrt{\frac{_{1-a}}{^2}}
\right\}
\]
for a real parameter $a \in (-1,1)$. This is an embedded surface in $S^3$ with mean curvature constant equal to
\[
H_a : =  \frac{2 \, a}{\sqrt{1-a^2}}.
\]
There are few other examples, and no general methods for the construction of CMC surfaces in $S^3$. The local operations involved in a gluing construction (such as forming connected sums using small bridging surfaces near a point of mutual tangency of two surfaces) have straightforward generalizations and gluing constructions can be carried out  in other ambient manifolds. However, the global aspects of the gluing construction are more complicated and restrictive when the examples one wishes to construct in this way are compact.

\paragraph*{Doubling the Clifford Torus.}
The construction that will be carried out in this paper produces new examples of embedded, higher-genus CMC surfaces of $S^3$, with small but non-zero mean curvature, by \emph{doubling} the unique minimal Clifford torus $C_0$ in the family of Clifford tori of $S^3$. This construction begins with the observation that the Clifford tori form a parallel foliation of a tubular neighbourhood of $C_0$ that is parameterized symmetrically on either side of $C_0$ by the mean curvature. The doubling construction consists in connecting together  two Clifford tori $C_a$ and $C_{-a}$ that have opposite, small mean curvature and are  symmetrically located on either side of $C_0$ at a each point of a sub-lattice of $C_0$. This gluing is performed using small bridging surfaces shaped like properly re-scaled catenoids  whose axis are perpendicular to the two initial tori and pass through the points of the sub-lattice. The resulting surface (henceforth called $\tilde S_a$) is topologically non-trivial, with genus $1+m_\Lambda$ where $m_\Lambda$ is the number of points of the sub-lattice.

Of course, $\tilde S_a$ does \emph{not} yet have constant mean curvature: its mean curvature is constant everywhere except near the gluing points.  One now hopes to find the desired CMC surface by \emph{perturbing} the surface $\tilde S_a$; this is done by expressing surfaces near to $\tilde S_a$ as normal graphs over $\tilde S_a$ and solving the partial differential equation that
determines when the nearby surface has constant mean curvature. The usual obstruction, well-known to those who have studied gluing constructions, arises at this point: in general, it turns out that this partial differential equation at the linearized level is  not bijective with a uniformly bounded inverse due to the existence of a kernel as well as a finite-dimensional \emph{approximate kernel} constituted of eigenfunctions associated to zero or small eigenvalues. However, this obstruction can be avoided if one chooses the gluing points in a sufficiently symmetric way and if one is careful enough in constructing the approximate solution.

To state the result precisely, let us parameterize the Clifford torus $C_a$ by
\[
\phi_a (\mu_1 , \mu_2) : = \left( \sqrt{\frac{_{{1+a}}}{^2}} \, e^{i\, \mu_1} ,  \sqrt{\frac{_{{1- a}}}{^2}} \, e^{i\, \mu_2} \right)
\]
where $(\mu_1, \mu_2) = {\mathbf R}^2$.  We will often identify ${\mathbf R}^4$ with ${\mathbf C}^2$. We consider $\Lambda$ a lattice of the plane $\Lambda : =  \{ n_1 \, \tau_1 + n_2 \, \tau_2 \, : \, n_1, n_2 \in {\mathbf Z}\}$ that is generated by $\tau_1 : = (\alpha_1, \beta_1)$ and $\tau_2 : = (\alpha_2, \beta_2)$
We assume that $\Lambda$ contains $2 \, \pi \, {\mathbf Z}^2$ so that $\phi_a (\Lambda)$ is a sublattice of $C_a$ and we denote by $m_\Lambda$ the number of points of $\phi_a (\Lambda)$ (which is also equal to the number of points of $\Lambda$ in $[0,2\pi)^2$).  

We denote by $G \subset SO(4)$, the group generated by
\[
\sigma_j (z_1, z_2) = (\me^{\mi \alpha_j}  \, z_1, \me^{\mi \beta_j} \, z_2)
\]
for $j=1,2$ and
\[
\rho (z_1, z_2) = (\bar z_1, \bar z_2).
\]
The theorem that this paper proves is the following.
\begin{thm}
Assume that $\Lambda$ is not contained in the curve $\{ (\mu_1, \mu_2 )  \, : \, \mu_1+ \mu_2 \equiv 0 \quad [2 \,\pi]\}$ or in the curve $\{ (\mu_1, \mu_2 )  \, : \, \mu_1 - \mu_2 \equiv 0 \quad [2 \, \pi] \}$. Then for all sufficiently small $a >0$, there exists a smooth, embedded surface $S_a \subset S^3$ with the following properties.
\begin{enumerate}
\item The surface $S_a$ is invariant under the action of $G$.

\item The surface $S_a$ has genus $m_\Lambda+1$.

\item The mean curvature of  $S_a$ is constant equal to $ H_a = \frac{2 \, a}{\sqrt{1-a^2}}$.

\item The surface $S_a$ is a small (normal) perturbation of the surface consisting of the connected sum of the two Clifford tori $C_a$ and  $C_{-a}$ at the points of $\phi_{a}(\Lambda)$ and  $\phi_{-a}(\Lambda)$ using small catenoids centered at the points of $\phi_0(\Lambda)$.

\item  As $a$ approaches zero, $S_a$ converges in ${\mathcal C}^\infty_{loc}$ topology to two copies of $C_0$ away from the points of $\phi_0(\Lambda)$.
\end{enumerate}
\label{th:1}
\end{thm}

Constant mean curvature surfaces of $\R^3$ analogous to the ones proposed above have been shown to exist by  Ritor\'e in \cite{ritore}, in which he constructs unbranched CMC immersions from a Riemann surface to $\R^3 / \Gamma$, where $\Gamma$ is a discrete group of translations of rank $g = 1, 2$ or $3$ generated by orthogonal vectors.  These immersions can often be lifted to periodic CMC immersions into $\R^3$.  In the rank $2$ case, these immersions are two-periodic and have a reflection symmetry with respect to the plane containing the generating vectors of the
translations. Furthermore, the fundamental domain of each immersed surface (i.e. the smallest domain which can be extended by periodicity to cover the entire surface) can be parametrized over a tubular neighbourhood around the ``$+$''-shaped one-dimensional variety formed by the union of two orthogonal orbits of the translations, and in some cases the upper and lower parts of the immersed surfaces are graphical over these tubular neighbourhoods. Ritor\'e does not use gluing techniques to construct these surfaces; rather he uses the Weierstra\ss-type representation of CMC surfaces. Moreover, the Karcher and Gro{\ss}e-Brauckmann examples, constructed by Schwartz reflection, duplicate some of Ritor\'e's examples and generate other examples that can also be considered analogous to those of the Main Theorem.

A further analogous construction has been achieved by Kapouleas and Yang in \cite{sdyang}, in which they construct a \emph{minimal} surface (i.e.~a surface with constant \emph{zero} mean curvature) by doubling the minimal Clifford torus in a similar manner as in Theorem~\ref{th:1}, except that they use an extremely large number of small bridging surfaces centered on a sub-lattice of very high order for the gluing.  (This result can be subsumed into the work of Pitts and Rubinstein in \cite{pittsrubinstein} for constructing equivariant minimal surfaces in space forms using symmetry and minimax methods.)  Kapouleas and Yang have discovered that it is possible to perturb the approximate surface to have exactly zero mean curvature when a particular relation between the number of catenoidal bridges and the size of the perturbation parameter is satisfied --- which can occur only for large numbers of catenoids and small perturbation parameter.  Thus in this framework, the obstructions to the perturbation can be avoided without the extra flexibility provided by allowing the mean curvature to vary, but the analysis is in this case much more delicate.

\paragraph*{Outline of the Proof.}
Theorem~\ref{th:1} will be proved in the following way.  One first expresses a small perturbation of the approximately CMC surface described above as a normal graph over $\tilde S_a$ whose graphing function $f$ belongs to a suitable Banach space.  Such a surface has the form $\mathrm{exp}(f N)(\tilde S_a)$ where $N$ is a smooth choice of unit normal vector field for $\tilde S_a$ and $\mathrm{exp}$ is the exponential map of the ambient $S^3$.  One then hopes to select a function $f_a$ which solves the partial differential equation $H \big( \mathrm{exp}(f N)(\tilde S_a) \big) = H_a$, where $H(\cdot)$ is the mean curvature function (with respect
to $N$), so that $\exp(f_a N)(\tilde S_a)$ is the desired CMC perturbation. One would accomplish this by applying the inverse function theorem to the non-linear partial differential operator
$\Phi_a(f) := H \big( \mathrm{exp}(f N)(\tilde S_a) \big) - H_a$ near $f=0$, which states that if the \emph{linearization} of $\Phi_a$ at $f=0$ is bijective with uniformly bounded inverse, then $\Phi_a$ can be inverted on a small neighbourhood of $\Phi_a(0)$.  Thus if $\Phi_a(0)$ is sufficiently small --- i.e.~that the mean curvature of $\tilde S_a$ deviates very little from $H_a$ --- then there exists $f_a$ so that $\Phi_a (f_a) = 0$.

Unfortunately, $\Dif \Phi_a (0)$ is not bijective with uniformly bounded inverse on an arbitrary Banach space and so the inverse function theorem does not apply in general, for two reasons.
First, the isometries of the ambient $S^3$ preserve mean curvature and thus all infinitesimal isometries are in the kernel of $\Dif \Phi_a (0)$. Second, when the surface consists of several
constituent pieces separated by small necks, as in the present case, then those motions of the surface corresponding to an infinitesimal isometry on one of the constituents and keeping the others fixed
(with transition in the neck regions), generate for  $\Dif \Phi_a (0)$ small eigenvalues tending to zero as $a \rightarrow 0$.  These two phenomena ensure that $\Dif \Phi_a (0)$ fails to be bounded below by a constant that does not tend to zero as $a \rightarrow 0$ on any Banach space that is not transverse to the kernel and approximate kernel of $\Dif \Phi_a (0)$.  The name \emph{Jacobi field} has been given to elements of the kernel of $\Dif \Phi_a (0)$ arising from geometric motions as described above.

If additional assumptions are made about the placement of the gluing points, then the  obstructions to controllable invertibility can be avoided.  Indeed, if the gluing points are located with sufficient symmetry and $\tilde S_a$ is deformed equivariantly (i.e.~deformations of $\tilde S_a$ are forced to preserve all symmetries), then the controllable invertibility of $\Dif \Phi_a (0)$ is contingent on whether the Jacobi fields --- both the global ones and those on the individual constituents of $\tilde S_a$ --- possess these additional symmetries or not. If, on each summand of $\tilde S_a$, there are no Jacobi fields possessing the symmetries, then the space of equivariant deformations of $\tilde S_a$ is transverse to the kernel and approximate kernel associated to small eigenvalues, and thus $\Dif \Phi_a (0)$ is controllably invertible.

\section{The Approximate Solution}

The purpose of this section of the paper is to construct the approximate solution $\tilde S_a$ and derive its relevant geometric properties.  This begins with a careful description of the building
blocks that will be assembled to construct $\tilde S_a$~: the Clifford tori in $S^3$ and the catenoid in $\R^3$. Since the proof of  Theorem~\ref{th:1} hinges on being able to rule out the existence of Jacobi fields and the approximate Jacobi fields of $\tilde S_a$, careful attention will be paid to understanding the Jacobi fields in each case.

\subsection{The Mean Curvature Operator and its Jacobi Fields}

The reader should be reminded of the linearized mean curvature operator of an arbitrary surface and of the origin of its `geometric' Jacobi fields.  Let $\Sigma$ be a closed hypersurface in a Riemannian manifold $X$ with mean curvature $H_\Sigma$, second fundamental form $b_\Sigma$ and unit normal vector field $N_\Sigma$. The linearization of the mean curvature operator on the space of normal graphs over $\Sigma$ is usually referred to as the Jacobi operator about $\Sigma$ and is given by
\[
{\mathcal L}_\Sigma : = \Delta_\Sigma + \big( \Vert B_\Sigma \Vert^2 + {\mathrm{Ric}} \, (N_\Sigma, N_\Sigma) \big)
\]
where $\Delta_\Sigma$ is the Laplace operator of $\Sigma$ and ${\mathrm{Ric}}$ is the Ricci tensor of $X$.  If $R_t$ is a one-parameter family of isometries of $X$ with deformation vector field $V = \left. \frac{\dif }{\dif t} \right|_{t=0} R_t$ then the function $\langle V, N \rangle$ is a solution of ${\mathcal L}_\Sigma \langle V, N \rangle =0$.  The function $\langle V, N \rangle$ is a \emph{Jacobi field} of $\Sigma$.

When $X=S^3$, the linearized mean curvature reads 
\[
{\mathcal L}_\Sigma  = \Delta_\Sigma + \big(\Vert B_\Sigma \Vert^2 + 2 \big) 
\] and the isometries of $S^3$ are simply the $SO(4)$-rotations of the ambient $\R^4$.  Thus there is at most a $6$-dimensional space of `geometric' Jacobi fields of $\Sigma$ arising from isometries. Since one expects in general that hypersurfaces with fixed constant mean curvature are isolated up to isometries, one expects no other Jacobi fields than these `geometric' ones.

\subsection{The Clifford Tori in \boldmath{$S^3$}}

The sphere $S^3$ contains a family of constant mean curvature surfaces known as the Clifford tori $C_a$ given by
\[
C_a : =  \left\{(z_1, z_2) \in \C^2 \, : \, |z_1| = \sqrt{\frac{_{1+a}}{^2}}, \quad |z_2| = \sqrt{\frac{_{1-a}}{^2}} \right\} ,
\]
for a real parameter $a \in (-1,1)$.  We now list the facts about the intrinsic and extrinsic geometry of the Clifford tori that will be needed in this paper.  In short, $C_a$ is \emph{flat} and has \emph{constant mean curvature}. To parameterize $C_a$ we define $\phi_a : S^1 \times S^1 \rightarrow S^3$ by
\[
\phi_a(\mu_1, \mu_2) : = \left(  \sqrt{\frac{_{1+a}}{^2}} \, e^{\mi \mu_1},  \sqrt{\frac{_{1- a}}{^2}} \, e^{\mi \mu_2} \right).
\]
Then with respect to the coordinates induced by $\phi_a$, the induced metric of $C_a$ is
\[
g_a : = \phi_a^\ast \left( ( \dif z_1 )^2 + ( \dif z_2)^2 \right) = \frac{_{1+a}}{^2} \, (\dif \mu_1)^2 + \frac{_{1-a}}{^2} \, (\dif \mu_2)^2
\]
so that $C_a$ is flat. If the unit normal vector field of $C_a$ is  chosen to be 
\[
N_a : = \sqrt{\frac{_{1-a}}{^2}}  \, e^{i\mu_1} \, \partial_{z_1} - \sqrt{\frac{_{1+a}}{^2}}
 \, e^{i\mu_2} \, \partial_{z_2} ,
\] 
then the second fundamental form of $C_a$ is
\[
B_a : = \frac{\sqrt{_{1-a^2}}}{^2} \, ( (\dif \mu_2)^2  - (\dif \mu_1)^2 )
\]
and the mean curvature of $C_a$ is the constant equal to 
\[
H_a : = \frac{_{2a}}{^{\sqrt{1-a^2}}}.
\] 
Note that the Clifford torus corresponding to $a = 0$ is the unique minimal submanifold in the family $C_a$.

\paragraph*{Jacobi fields.}  The linearization of the mean curvature operator on the space of normal graphs over $C_a$ is given by 
\[ 
\mathcal{L}_{C_a}  : = \frac{_2}{^{1+a}} \, \partial_{\mu_1}^2 + \frac{_2}{^{1-a}} \,  \partial_{\mu_2}^2  + \frac{_4}{^{1-a^2}} 
\]
Two results concerning the Jacobi fields of the operator $\mathcal{L}_{C_0}$ will be needed in the sequel.

\begin{lemma}
    \mylabel{lemma:jacobi}
    The Jacobi fields of $C_0$ are linear combinations of functions of the form $(\mu_1 ,\mu_2) \longmapsto \cos (\mu_1 \pm \mu_2)$ and $(\mu_1 ,\mu_2) \longmapsto \sin (\mu_1 \pm \mu_2)$.
\end{lemma}

\begin{proof}
This result follows at once from Fourier analysis in both the $\mu_1$ and $\mu_2$ variables.  Observe that any of these these Jacobi fields is  associated to a one parameter family of isometries of $S^3$.
\end{proof}

We now consider $\Lambda$ a sub-lattice of the plane $\Lambda : =  \{ n_1 \, \tau_1 + n_2 \, \tau_2 \, : \, n_1, n_2 \in {\mathbf Z}\}$ generated by $\tau_1 : = (\alpha_1, \beta_1)$ and $\tau_2 : = (\alpha_2, \beta_2)$.  We assume that $\Lambda$ contains $2 \, \pi \, {\mathbf Z}^2$.

\begin{cor}
    \label{cor:nojacobi}
Assume that $\Lambda$ is not contained in $$\{ (\mu_1,  \mu_2 ) \, : \, \mu_1+ \mu_2 \equiv 0 \: [2\pi]\} \qquad \mbox{or} \qquad  \{ (\mu_1, \mu_2 ) \, : \, \mu_1 - \mu_2 \equiv 0 \: [2\pi]\}\, .$$ Then there is no nontrivial function $f$ that solves the equation  ${\mathcal L}_{C_0} \, f =0$ and that satisfies
\[
f( \mu+ \tau_j ) = f(\mu) \qquad \mbox{and} \qquad f(-\mu) = f(\mu)
\] 
for all $\mu := (\mu_1, \mu_2) \in {\mathbf R}^2$.
\end{cor}
\begin{proof}
Using the previous lemma we know that $f$ is a linear combination of the functions $(\mu_1 ,\mu_2) \longmapsto \cos (\mu_1 \pm \mu_2)$ and $(\mu_1 ,\mu_2) \longmapsto \sin (\mu_1 \pm \mu_2)$. But since we assume that $f(-\mu) = f(\mu)$ for all $\mu$, then we conclude that $f$ is a linear combination of $(\mu_1 ,\mu_2) \longmapsto \cos (\mu_1 \pm \mu_2)$.

We write
\[
f(\mu_1 ,\mu_2) =  a^+ \, \cos (\mu_1 + \mu_2) + a^- \cos(\mu_1 - \mu_2)
\]
Let us assume that $a^+\neq 0$ and $a^- \neq 0$. Observe that $(0,0)$ is a critical point of $f$  and the second variation of $f$ at $(0,0)$ is given by
\[
D^2f_{(0,0)} (\mu_1, \mu_2) =  - \frac{_1}{^2} \, \left( a^+ \, (\mu_1 \pm \mu_2)^2 + a^- (\mu_1 - \mu_2)^2\right) .
\]
The only other critical points of $f$ are of the form $( (n_1+n_2) \, \frac{\pi}{2} , ( n_1-n_2 ) \, \frac{\pi}{2})$, where $n_1, n_2 \in {\mathbf Z}$. However among them, the only critical point that has the same second variation are the points $( (n_1+n_2) \,  \pi  ,(n_1-n_2 ) \,  \pi )$, where $n_1, n_2 \in {\mathbf Z}$.  Now, if we use the fact that $f$ satisfies $f( \mu + \tau_j ) = f(\mu)$, for $j=1,2$, we find that the lattice $\Lambda$ has to be included in $\pi \, {\mathbf Z}^2$. This contradicts the hypothesis that the lattice $\Lambda$ is not included in $\{(\mu_1,  \mu_2 ) \, : \, \mu_1 + \mu_2 \equiv 0 \: [2\pi]\}$.

We now assume that $a^+ =1$ and $a^-=0$. Then, the set of points where $f=1$ is precisely equal to $\{ (\mu_1 , \mu_2 ) \, : \, \mu_1 + \mu_2 \equiv 0 \: [2\pi]\}$. Again, if we use the fact that $f$ satisfies $f( \mu + \tau_j ) = f(\mu)$, for $j=1,2$,  this contradicts the hypothesis that the lattice $\Lambda$ is not included in $\{ (\mu_1 , \mu_2 ) \, : \, \mu_1 + \mu_2 \equiv 0 \: [2\pi]\}$. Finally assume that $a^+=0$ and $a^-=1$. Then the set of points where $f=1$ is precisely equal to $\{ (\mu_1 , \mu_2 ) \, : \, \mu_1 - \mu_2 \equiv 0 \: [2\pi]\}$. Again this contradicts the hypothesis that the lattice $\Lambda$ is not included in $\{ (\mu_1 , \mu_2 ) \, : \, \mu_1 - \mu_2 \equiv 0 \: [2\pi]\}$. This implies that $a^+ = a^- =0$ and hence $f \equiv 0$.
\end{proof}

\subsection{The Catenoid in  \boldmath{${\mathbf R}^3$}}
\label{sec:catenoid}

The catenoid $K$ is the unique, complete, two-ended and cylindrically  symmetric embedded minimal surface in $\R^3$.  Re-scalings of $K$ are  also minimal surfaces. The \emph{$\eps$-scaled catenoid} is the image of  ${\mathbf R}\times S^1$ by 
\[
\psi_\eps (s, \theta) : =  \left( \eps \, \cosh s \, \cos \theta, \eps \, \cosh s \, \sin \theta, \eps \, s  \right) .
\]
In this parametrization  the induced metric of $\eps K$ is 
\[
g_{\eps K} := \psi_\eps^* \,(\dif x^2 + dy^2 + dz^2) = \eps^2 \, \cosh^2 s\, \big( \dif s^2 + \dif \theta^2  \big) \, .
\]
If the unit normal vector field of $\eps K$ is chosen to be
\[
N_{\eps K} : = \frac{_1}{^{\cosh s}} \, \left(   - \cos \theta \, \partial_x - \sin \theta \, \partial_y + \sinh s \, \partial_z \right) \, ,
\] 
then the second fundamental form of $\eps K$ is given by
\[
B_{\eps K} = \eps \, ( d\theta^2 - ds^2) ,
\]
and the mean curvature of $\eps K$ vanishes.  

\paragraph*{Jacobi fields.} The linearization of the mean curvature operator on the space of normal graphs over $\eps K$ is here given by 
\[
\mathcal{L}_{\eps K}  =  \frac{_1}{^{\eps^2 \, \cosh^2 s}} \, \left( \partial_s^2 +\partial_\theta^2 + \frac{2}{\cosh^2 s} \right)  .
\]
The Jacobi fields of $\eps K$ are the solutions of the equation $\mathcal{L}_{\eps K} \, u  = 0$.   Some Jacobi fields can be explicitly computed since they are associated to one parameter families of minimal surfaces to which $\eps K$ belongs. For example, the translations 
and rotations in ${\mathbf R}^3$ generate five linearly independent Jacobi fields and dilation generates a sixth Jacobi field. 
\begin{lemma}
    \mylabel{lemma:catenoidjac}
Assume that $\delta <2$ is fixed. The subspace of Jacobi fields of $\eps K$ that are 
bounded by a constant times $(\cosh s)^\delta $  and are invariant under the action of the symmetry $(s , \theta) \mapsto (s,  \theta +\pi )$ is two-dimensional spanned by the functions $f^{(1)}_0 (s, \theta) : = \tanh s$ and $f^{(2)}_0(s, \theta) := s \, \tanh s - 1$.
\end{lemma}

\begin{proof}
We decompose the Jacobi field $f$ in Fourier series and write
\[ 
f(s, \theta) = \sum_{j=0}^\infty (a_j (s) \, \cos (j\theta)+ b_j(s) \, \sin (j\, \theta)) .
\]
Since we assume that $f(s , \theta) =f (s, \theta +\pi)$, we get  
\[
f(s, \theta) = \sum_{j=0}^\infty (a_{2j} (s) \, \cos (2j \, \theta)+ b_{2j}(s) \, \sin (2j\, \theta)) .
\]

Now, when $j \geq 1$ the functions $a_{2j}$ and $b_{2j}$ are solutions of the ordinary differential equation
\[
\left(\partial_s^2  -(2j)^2 + \frac{_2}{^{\cosh^2 s}}  \right) \, a_{2j} = \left(\partial_s^2  -(2j)^2 + \frac{_2}{^{\cosh^2 s}} \right) \, b_{2j}=0 .
\]
The study of the possible asymptotics of $a_{2j}$ and $b_{2j}$ at $\pm \infty$ show that either these functions blow up like $(\cosh s)^{2j}$ or decay like $(\cosh s)^{-2j}$. However, since they are bounded by a constant times $(\cosh s )^\delta$ for some $\delta <2$, they have to decay exponentially at $\pm \infty$. Then the maximum  principle implies that $a_{2j}= b_{2j}=0$  since, for $j \geq 1$, the potential of this ordinary differential equations is negative.

When $j=0$, two independent solutions of the homogeneous problem 
\[
\left(\partial_s^2  + \frac{_2}{^{\cosh^2 s}}  \right) \, f_0 = 0 .
\]
are known. The Jacobi fields associated to vertical translations of the ambient $\R^3$ is given explicitly as
$$
       f_0^{(1)}(s, \theta) := \big\langle {\partial_z}, N_{\eps K} \big\rangle = \tanh s
$$
and there is one other Jacobi field coming from the dilations of $\eps K$ and is given explicitly by
$$ f_0^{(2)}(s, \theta) := \big\langle \psi_{\eps K}, N_{\eps K} \big\rangle = \eps \, \left(  s \, \tanh s - 1   \right)\, .$$ Obviously the function $f_0$ is a linear combination of these two functions since it is a solution  of a linear second order ordinary differential equation. This completes the proof of the result.
\end{proof}

\subsection{Toroidal Coordinates for $S^3$}

The approximate solution of the torus doubling construction will be formed by taking two Clifford  tori of the form $C_{\pm a}$ symmetrically placed on either side of the unique minimal Clifford  torus $C_0$, and connecting them at a sub-lattice of points by bridging surfaces consisting of small pieces of re-scaled catenoids, embedded into $S^3$.  In order to perform  this construction with as much precision as possible, it is most convenient to use canonical coordinates for a neighbourhood of $C_0$ in $S^3$ which are well-adapted to the family of  Clifford tori and that can be used to embed the catenoids with the least amount of distortion.

We set $S^3_* : =  S^3 \setminus \left( S^1\times  \{0\} \cup \{0\} \times S^1\right)$. Let $I = (-\frac{\pi}{2\sqrt 2},\frac{\pi}{2\sqrt 2})$.  The \emph{toroidal coordinates for $S^3$} are given by the inverse of the mapping $\Xi : S^1 \times S^1 \times I \rightarrow S^3_*$   (written as a subset of ${\mathbf C}^2$) which is defined by
    $$\Xi (\mu_1, \mu_2,  t) : = \left( \sqrt{ \frac{_{1+\sin (\sqrt 2 t)}}{^2}}  \, e^{\mi  \mu_1}, 
    \sqrt{ \frac{_{1-\sin (\sqrt 2 t)}}{^2}}   \, e^{\mi  \mu_2} \right) \, .$$

\begin{lemma}
    The standard metric on $S^3$ expressed in toroidal coordinates is
 \[
 \Xi^\ast g_{S^3} =  \frac{_1}{^2}\, \dif t^2 +  \frac{_{1+ \sin (\sqrt 2 t)}}{^2} \,  ( \dif \mu_1)^2 +  \frac{_{1-\sin (\sqrt 2 t)}}{^2}   \,( \dif \mu_2)^2\, .
 \]
\end{lemma}

Close to $C_0$, $t$ is small and the metric $\Xi^\ast g_{S^3}$ can be considered a perturbation of the flat metric
 \[
\mathring g : = \frac{_1}{^2} ( \dif t^2 + (\dif \mu^1)^2 +  (\dif \mu^2)^2).
\]   
This observation allows many (though not all) of the forthcoming estimates to be computed with respect to $\mathring g$, which will be simpler.

\subsection{Construction of the Approximate Solution}
\label{sec:construction}

We consider $\Lambda$ a lattice of the plane $\Lambda : =  \{ n_1 \, \tau_1 + n_2 \, \tau_2 \, : \, n_1, n_2 \in {\mathbf Z}\}$ that is generated by $\tau_1 : = (\alpha_1, \beta_1)$ and $\tau_2 : = (\alpha_2, \beta_2)$.  We assume that $\Lambda$ contains $2 \, \pi \, {\mathbf Z}^2$, and further,  that the lattice $\Lambda$ is not contained in $\{ (\mu_1,  \mu_2 ) \, : \, \mu_1+ \mu_2 \equiv 0 \: [2\pi]\}$ or $\{ (\mu_1, \mu_2 ) \, : \, \mu_1 - \mu_2 \equiv 0 \: [2\pi]\}$.  According to the result of Corollary~\ref{cor:nojacobi} there is no non-trivial solution of ${\mathcal L}_{C_0} \, u=0$ that satisfies
\begin{equation}
u (\mu + \tau_j)=u (\mu) \quad \mbox{for $j=1,2$} \qquad \mbox{and}  \qquad u (-\mu) = u (\mu) .
\label{eq:eqi}
\end{equation}
This means that, restricted to the set of functions satisfying (\ref{eq:eqi}) the operator ${\mathcal L}_{C_0} :  { C}^{2,\alpha} (C_0) \longrightarrow { C}^{0, \alpha}(C_0)$ is an isomorphism.

The construction of the approximate solution begins with the definition of a function $\Gamma_\Lambda$ that is the unique solution of the equation
\[
{\mathcal L}_{C_0} \, \Gamma_\Lambda = - \pi \, \sum_{p \in \phi_0(\Lambda )} \delta_p
\]
satisfying (\ref{eq:eqi}). Here, $\delta_p$ is the Dirac $\delta$-mass at the point $p \in C_0$. With slight abuse of notations, we will write $\Gamma_\Lambda (\mu)$ instead of $\Gamma_\Lambda \circ \phi_0 (\mu)$. Observe that, near $\mu= (0,0)$ the function $\Gamma_\Lambda$ can be expanded as 
\[
\Gamma_\Lambda  (\mu) = - \log |\mu| + \gamma_\Lambda + {\mathcal O}(|\mu|^2 \, \log 1/|\mu|) 
\]
where $\gamma_\Lambda \in {\mathbf R}$ is a constant that depends on the lattice $\Lambda$. Here we have implicitly used the fact that $ \Gamma_\Lambda  (-\mu) = \Gamma_\Lambda  (\mu)$ so that the usual error term ${\mathcal O}(|\mu|) $ that appears in the Taylor expansion of $\Gamma_\Lambda$ can be improved into ${\mathcal O}(|\mu|^2\,  \log 1/|\mu| ) $.

Next, we consider the Clifford tori  $C_{\pm a}$ which are the images of $S^1\times S^1 \times \{\pm t_a\}$ under the toroidal coordinate mapping $\Xi$.  Here $t_a$ and $a$ are related by
\[
\sin (\sqrt 2 \, t_a ) = a \, .
\]
These Clifford tori will be perturbed using a proper multiple of the function $\Gamma_\Lambda$. To this end, we define the parameter $\eps_a >0$, for $a$ close enough to $0$, to be the unique positive solution of 
\[
t_a =  - \eps_a \, \log \eps_a + \eps_a \, ( \gamma_\Lambda +  \log 2 ) \, .
\]
We define ${\mathcal C}_a^+$ using the toroidal coordinates to be the image of ${\mathbf R}^2 \setminus \bigcup_{\mu \in \Lambda} B_{\sqrt {\eps_a}}(\mu)$ under the mapping 
\[
\mu \longmapsto \Xi(\mu, t_a - \eps_a \, \Gamma_\Lambda (\mu) )
\]
and also ${\mathcal C}_{a}^-$ to be the image of ${\mathbf R}^2 \setminus \cup_{\mu \in \Lambda} B_{\sqrt {\eps_a}}(\mu)$ under the mapping 
\[
\mu \longmapsto \Xi (\mu, - t_a + \eps_a  \, \Gamma_\Lambda (\mu) ) \, .
\]
This produces two surfaces that are close to $\phi_0 ({\mathbf R}^2 \setminus \cup_{\mu \in \Lambda} B_{\sqrt{\eps_a}}(\mu))$, and with $m_\lambda$ boundaries. We also consider the re-scaled catenoid $\eps_a \, K$ that we insert in $S^3$ as follows. Consider the image of $\{(s, \theta) \in {\mathbf R}\times S^1 \, : \, \eps_a \, \cosh s \leq \sqrt \eps_a\}$ by
\[
(s, \theta) \longmapsto \Xi (\eps_a \, \cosh s \, \cos \theta, \eps_a \, \cosh s \, \sin \theta, \eps_a s)
\]
as well as the images of this surface with boundaries by the action of the elements of the group $G$. This produces $m_\Lambda$ surfaces with boundaries the union of which will be denoted by ${\mathcal N}_a$. The union of these pieces of surfaces does not produce a smooth surface but using cut-off functions we can interpolate smoothly between the different summands ${\mathcal C}_a^\pm$ and ${\mathcal N}_a$ to obtain a smooth surface that will be denoted by $\tilde S_a$. 

Because of the invariance under the action of $G$ it is enough to explain how to connect the pieces of surfaces close to the point $\phi_0 (0,0)$. For example, near $\mu=(0,0)$,  the graph of $\mu \longmapsto  t_a - \eps_a \, \Gamma_\Lambda (\mu) $ can be expanded as 
\[
 t_a - \eps_a \, \Gamma_\Lambda (\mu) =  t_a + \eps_a \, \left( \log  |\mu|  -  \gamma_\Lambda + {\mathcal O} (  |\mu|^2\, \log  1/|\mu|) \right).
\]
While, changing variables, $|\mu |= \eps_a \, \cosh s $, with $s >0$, we find with little work that 
\[
\eps_a \, s (\mu) =  \eps_a \, \left( \log |\mu| +  \log (2/\eps_a) + {\mathcal O} (\eps_a^{2}\, |\mu|^{-2}) \right) = t_a +  \eps_a \, \left( \log |\mu|Ê- \gamma_\Lambda + {\mathcal O} (\eps_a^{2}\, |\mu|^{-2}) \right)\, .
\]
To obtain a smooth surface it is enough to interpolate between the  between these two functions inside an annulus of radii $2 \, \sqrt \eps_a$ and $\sqrt \eps_a /2$.  For example, to interpolate smoothly between the graph of $t_a - \eps_a \, \Gamma_\Lambda$ and the graph of $\eps_a \, s$ we define the function 
\[
T_{a} (\mu ) : =  \eta (\mu / \sqrt \eps_a) \, (t_a - \eps_a \, \Gamma_\Lambda (\mu)) + (1- \eta (\mu / \sqrt \eps_a) \, \eps_a \, s(\mu)
\]
where $\eta$ is a cut-off function identically equal to $0$ in $B_{1/2}(0)$ and identically equal to $1$ in ${\mathbf R}^2 \setminus B_2(0)$.  A similar analysis can be performed for the lower end of the re-scaled catenoid. The final step in the assembly of the different summands of the approximate solution is to extend the above construction so that the resulting surface is invariant under the action of the elements of $G$. We will denote the transition regions by ${\mathcal T}_a$ corresponding to the image of $\bar B_{2 \sqrt{\eps_a}} (0) \setminus B_{\sqrt{\eps_a}} (0)$ by $\mu \longmapsto \Xi(\mu, T_{a}(\mu))$. 

This recipe produces a surface $\tilde S_{a}$ that  is a smooth, embedded submanifold of $S^3$.  It is equal to the connected sum of $\eps_a$-re-scaled catenoids centered at the points of $\phi_ 0 (\Lambda)$ and small perturbations of the Clifford tori $C_{\pm a}$.  Recall that, by construction, these tori have mean curvature equal to $H_a := \frac{2a}{\sqrt{1-a^2}}$.  Finally, when $a$ approaches zero, then $\tilde S_a$ approaches two copies of the unique minimal Clifford torus, punctured at the sub-lattice of points $\phi_0(\Lambda)$.

The symmetries of the approximate solution $\tilde S_a$ constructed above will play a crucial role in the forthcoming analysis, where only deformations preserving these symmetries will be considered as valid.  This will have the effect of eliminating the kernel and approximate kernel of the Jacobi operator, which is the fundamental obstruction to the invertibility of the linearized deformation operator.  Thus one should observe that the approximate solution $\tilde S_a$ is by construction symmetric with respect to the rotations  $\sigma_{j}$ defining the sub-lattice as well as the symmetry $\rho$. Recall that these rotations are given by 
\[
\sigma_j (z_1,z_2) =  ( e^{i\alpha_j} \, z_1, e^{i\beta_j} \, z_2) ,
\]
while in toroidal coordinates of $S^3$ they are given by 
\[
\sigma_{j}\circ \Xi (\mu_1, \mu_2, t) =   \Xi (\mu_1+\alpha_j, \mu_2+\beta_j, t) \, .
\]
The symmetry $\rho$ is given by 
\[
\rho(z_1, z_2) = (\bar z_1, \bar z_2) ,
\]
while in toroidal coordinates of $S^3$ it is given by 
\[
\rho \circ \Xi  (\mu_1, \mu_2, t) = \Xi (-\mu_1 , -\mu_2, t) \, .
\]  
Since all of these symmetries preserve both $\tilde S_a$ and the ambient metric of $S^3$, these symmetries also preserve the mean curvature of $\tilde S_a$.

\subsection{Estimates for the Approximate Solution}

The remaining task for this section of the paper is to estimate all relevant geometric quantities on $\tilde S_a$ in terms of the parameter $a$.  Note that it is generally sufficient to estimate with respect to the geometry induced on $\tilde S_a$ by the flat metric $\mathring g$ since the geometry induced by the toroidal coordinates is a small perturbation of the Euclidean metric.  Also, by symmetry it is sufficient to estimate only at the sub-lattice point corresponding to $\mu = (0,0)$.  The most important of the estimates of $\papprox$ is the pointwise ${C}^{1}$ estimate of its mean curvature.
\begin{prop}
    \mylabel{prop:estH}
The mean curvature of the approximate solution $\papprox$ satisfies the following estimates.   If $\Xi (\mu, t ) \in \mathcal C_a^\pm$ then
\begin{subequations}
    \mylabel{eqn:meancurvest}
    \begin{equation}
        \mylabel{eqn:mcone}
        \vert H (\tilde S_a)  - H_a \vert + |\mu | \,  \Vert \nabla H  (\tilde S_a) \Vert  \leq   C \, ( \eps^2_a \, (\log 1/\eps_a) \, | \mu |^{-2}) \, ,
    \end{equation}
   if $\Xi (\mu, t) \in \mathcal T_a$ then
    \begin{equation}
        \mylabel{eqn:mctwo}
        \vert H (\tilde S_a)  - H_a  \vert + |\mu | \, \Vert \nabla H (\tilde S_a)  \Vert \leq  C \, \eps_a \, ( \log 1/\eps_a ) \, ,
    \end{equation}
and if $\Xi (\eps_a \cosh s \cos \theta, \eps_a \cosh s \sin \theta , \eps_a \, s) \in \mathcal N_a$ then
    \begin{equation}
        \mylabel{eqn:mcthree}
        \vert H (\tilde S_a)  - H_a  \vert + \eps_a \, \cosh s \, \Vert \nabla H (\tilde S_a) \Vert \leq  C \, (\log 1/\eps_a ) \, (\cosh s)^{-2} \, ,
    \end{equation}
\end{subequations}
where $C >0$ is a constant independent of $a$, provided $a$ is sufficiently small.
\end{prop}

\begin{proof}
The calculations can all be done using the toroidal coordinates, for which the ambient metric will be denoted $g$.  Let us use the abbreviations $\nabla_i = \nabla_{\partial_{\mu_i}}$ and $\nabla_t = \nabla_{\partial_t}$ and also use a comma to denote partial differentiation, such as $u_{,i} = \partial_{\mu_i} u$.    

To estimate the mean curvature in the regions $\mathcal  T_a$ and $\mathcal C_a^\pm$, we first compute the mean curvature of the graph of a function $\sqrt{2} u : S^1 \times S^1 \rightarrow \R$ parametrized by $\mu \longmapsto (\mu, \sqrt 2 \, u(\mu))$.  The tangent vectors of this surface are given by $T_j : =  \partial_{\mu_j} + \sqrt 2  u_{,j} \, \partial_t$ and it is easy to check that the induced metric is given by
\[
\bar g = \left( \frac{_1}{^2} ( 1+ \sin (2u)) + (u_{,1})^2 \right) \, d\mu_1^2
+ \left( \frac{_1}{^2} ( 1 - \sin (2u)) +( u_{,2})^2 \right) \, d\mu_2^2
+ 2 \, u_{,1} u_{,2} \, d\mu_1\, d\mu_2 \, .
\]
The normal vector field $N$ can be written as $N : = \bar N / ||Ê\bar N||$ where  $\bar N : =  \sqrt 2 \, \partial_t - a^j \, T_j$ and the coefficients $a^j$ are determined so that $\bar N$ is normal to the surface. One finds the explicit expressions $a^j =  g^{jk} u_{,k}$ and $||Ê\bar N||^2 =  1+ g^{jk} u_{,j} u_{,k}$.  We now compute
\begin{align*}
g(\nabla_{T_i} T_j,  \bar N) &=   \sqrt 2 g (\nabla_{T_i} T_j , \partial_t) - a^k g(\nabla_{T_i} T_j,  T_k)\\
	& =  \sqrt 2 g (\nabla_{T_i} T_j , \partial_t) -  \bar \Gamma^k_{ij} u_{,k}  
\end{align*}
where the $\bar \Gamma$-terms are the Christoffel symbols of the induced metric $\bar g$.   We evaluate the first term as
\begin{align*}
	g (\nabla_{T_i} T_j , \partial_t) &= g \big( (\nabla_i + \sqrt{2} u_{,i} \nabla_t) (\partial_j + \sqrt{2} u_{,j} \partial_t), \partial_t \big) \\
	&= \Gamma_{ijt} + \tfrac{1}{\sqrt{2}} u_{,ij} + \sqrt{2}u_{,j}  \Gamma_{itt} + \sqrt{2} u_{,i} \Gamma_{tjt} + 2 u_{,i} u_{,j} \Gamma_{ttt} 
\end{align*}
where the $\Gamma$-terms are the Christoffel symbols of the ambient metric $g$.  Then  $\Gamma_{11t} = - \tfrac{1}{2\sqrt{2}} \cos (\sqrt{2} t)$ and $\Gamma_{22t} =  \tfrac{1}{2\sqrt{2}} \cos (\sqrt{2} t)$, whereas all other $\Gamma$-terms vanish.  We thus obtain the second fundamental form
\[
 \, \bar B =  \frac{_1}{^{||\bar N||}} \left(\tfrac{1}{2} \cos (2u) \, ( d\mu_2^2 - d\mu_1^2) +  u_{,ij}  \, d\mu_i\, d\mu_j - \bar \Gamma^k_{ij} u_{,k}   \, d\mu_i\, d\mu_j\right) \, .
\] 
Finally, by taking the trace of $\bar g^{-1} \bar B$ we get the mean curvature
\begin{equation}
	\label{eqn:secondff}
	\begin{aligned}
		H &= \frac{1}{2 \det(\bar g)  ||\bar N||} \Big[  u_{,11} \big( 1 - \sin(2u) + 2 (u_{,2})^2 \big) - 4 u_{,1} u_{,2} u_{,12} + u_{,22} \big( 1 + \sin(2u) + 2 (u_{,1})^2 \big)  \\
		&\qquad \qquad \quad \qquad +    \sin(2 u) \cos(2 u)  + \cos(2u) \big( (u_{,1})^2 - (u_{,2})^2 \big) \Big] + \frac{1}{||\bar N||} \bar g^{ij} \bar \Gamma_{ij}^k u_{,k} 
	\end{aligned}
\end{equation}
where $\det(\bar g) =\frac{1}{4} \big( \cos^2(2u) + 2(1-\sin(2u)) (u_{,1})^2 + 2 (1+\sin(2u) )(u_{,2})^2\big)$.

Consider now the region ${\mathcal C}^+_a$ in which case $u=u_a$ with $\sqrt 2 \, u_a : = t_a - \eps_a \, \Gamma_\Lambda$ and $|\mu | \geq \frac{1}{2} \, \sqrt{\eps_a}$.   First, we can check that $||Ê\bar N||  =  1 + {\mathcal O}(\eps^2_a |\mu|^{-2})$ and $\bar g_{ij} = \delta_{ij} +  {\mathcal O}(\eps^2_a |\mu|^{-2})$, so that $\bar g^{ij} \bar \Gamma^k_{ij} u_{,k} = {\mathcal O}(\eps^2_a |\mu|^{-2})$ as well, since this quantity consists of terms that are of the form $u_{,i} u_{,j}$ or  $u_{,i} u_{,j} u_{,kl}$ multiplied by coefficients of $\bar g^{-1}$.  Now, using the fact that $D u_a = \mathcal O (\eps_a |\mu|^{-1})$ and $D^2 u_a = \mathcal O(\eps_a |\mu|^{-2})$ we find 
$$H = 2 \, \tan (2u_a)  + 2\, \Delta u_a +{\mathcal O}(\eps^2_a |\mu|^{-2})   \, .$$
Finally, we use the fact that $(\Delta +2) \, \Gamma_\Lambda  =0$ away from the points of $\Lambda$ and also the formula $H_a = 2 \tan (\sqrt{2} \, t_a)$ and the estimate $u_a ={\mathcal O}(\eps_a \, (\log 1/\eps_a))$  to conclude $H =  H_a +{\mathcal O}(\eps^2_a \, (\log 1/\eps_a) \, |\mu|^{-2})$ holds in $\mathcal C^-_a$.  The corresponding estimates in ${\mathcal C}^-_a$ and also in ${\mathcal T}_a$ are obtained using similar computations, as well as the estimate of the derivative of $H$.

It remains to compute the mean curvature of the neck region ${\mathcal N}_a$. Since the center of the neck region is not a graph over the level sets of constant $t$, the previous calculation does not help us.  Thus we compute directly the mean curvature of the surface parametrized by  
\[
(s, \theta) \longmapsto (\eps_a \, \cosh s\, \cos \theta , \eps_a \, \cosh s\, \sin \theta, \eps_a  s)
\]
for $s$ such that $\eps_a \cosh s \leq \frac{1}{2}Ê\, \sqrt{\eps_a}$.  The tangent vectors of this surface are given by 
\[
T_\theta := \eps_a \, \cosh s  \,  (\cos \theta \, \partial_{\mu_2} -  \sin \theta \, \partial_{\mu_1}) \qquad \mbox{and} \qquad T_s : =  \eps_a \, \sinh s  (\cos \theta \, \partial_{\mu_1} + \sin \theta \, \partial_{\mu_2} ) + \eps_a \, \partial_t .
\] 
The induced metric is then given by
\begin{align*}
	\bar g &=  \frac{_1}{^2}\, \eps^2_a \, \cosh^2 s  \, (ds^2 +d \theta^2) +   \frac{_1}{^2}  \, \eps^2_a \,  \sin (\sqrt 2 \eps_a s) \, \cos (2\theta)  \, (  \sinh^2 s \, ds^2 -  \cosh^2 s \, d\theta^2) \\
	&\qquad  -   \eps_a^2 \, \cosh s\, \sinh s \, \sin(\sqrt 2 \, \eps_a\, s) \, \sin (2\theta) \,  ds \, d\theta   \, .
\end{align*}
The normal vector field $N$ can be written as $N : = \bar N / ||Ê\bar N||$ with $\bar N : = N_0 - a^s \, T_s - a^\theta \, T_\theta$ where 
$$N_0 : =   \tfrac{1}{\cosh s} \, \left( -\cos \theta \, \partial_{\mu_1} - \sin \theta \, \partial_{\mu_2} + \sinh s \, \partial_tÊ\right) $$
and where the coefficients $a^s$, $a^\theta$ are determined so that $\bar N$ is normal to the surface.  One finds the estimates
\begin{align*}
	a^s &=  {\mathcal O} \left(\tfrac{\log (1/\eps_a)}{ \cosh^2 s} \right)\\
	a^\theta &=  {\mathcal O} \left(\tfrac{\log (1/\eps_a)}{ \cosh^2 s} \right) \\
	||Ê\bar N||^2 &=  \tfrac{1}{2} +  {\mathcal O} \left(\tfrac{\eps_a \, (\log (1/\eps_a))}{ \cosh^2 s} \right) \, .
\end{align*}

We now compute the coefficients of the second fundamental form.  Our calculations are simplified by considering the local parametrization 
\[
\psi : (s, \theta, x) \longmapsto \left( ( \eps_a \, \cosh s - \tfrac{x}{\cosh s} )\, \cos \theta , (\eps_a \, \cosh s - \tfrac{x}{\cosh s} ) \, \sin \theta, \eps_a  s + x \, \tanh s \right)
\]
so that $N_0  = \psi_\ast \, \partial_x$.  Hence,
\begin{align*}
	- 2 ||\bar N || B(T_i, T_j) &= g(\nabla_{T_i} \bar N, T_j) + g(\nabla_{T_j} \bar N, T_i) \\
	 &= g(\nabla_{T_i} N_0, T_j) + g(\nabla_{T_j} N_0, T_i) + a^k_{,i} \bar g_{jk} + a^k_{,j} \bar g_{ik} + a^k \bar g_{ij,k}\\
	 &=  \partial_x g(T_i, T_j) |_{x=0} +  a^k_{,i} \bar g_{jk} + a^k_{,j} \bar g_{ik} + a^k \bar g_{ij,k} 
\end{align*}
where $i, j, k$ can be $s$ or $\theta$.  It is easy to check that $a^k_{,i} \bar g_{jk} + a^k_{,j} \bar g_{ik} + a^k \bar g_{ij,k}  = {\mathcal O} ( \eps_a^2 \, (\log (1/\eps_a))$.  Finally,
\[
\partial_x g(T_i, T_j) |_{x=0} = 
\begin{cases}
	 \eps_a + \frac{1}{\sqrt 2} \, \eps_a^2 \, \cosh^2 s \, \tanh s \, \cos (\sqrt 2 \, \eps_a s) \cos (2 \theta) \\[1ex] 
	\qquad \qquad  + \eps_a \, \tanh^2 s \, \sin (\sqrt{2} \eps_a s) \, \cos ( 2\theta) \quad \mbox{when $i=j=s$},  \\[1ex]
	- \eps_a - \frac{1}{\sqrt 2} \, \eps_a^2 \, \cosh^2 s \, \tanh s \, \cos (\sqrt 2 \, \eps_a s) \cos (2 \theta) \\[1ex]
	\qquad \qquad + \eps_a \, \sin (\sqrt{2} \eps_a s) \, \cos ( 2\theta) \quad \mbox{when $i=j=\theta$},\\[1ex]
 	- \frac{1}{2} \, \eps_a^2 \, \cosh s\, \sinh s \, \tanh s \, \sin (2 \theta)  \, \cos ( \sqrt 2 \, \eps_a\, s)  \\[1ex]
	\qquad \qquad \mbox{ when $i\neq j$} \, .
\end{cases}
\]
For example, 
\[
g(T_s, T_\theta) = \frac{_1}{^2} \, (\eps_a+ \frac{_x}{^{\cosh s}})^2 + \frac{_1}{^2} \, ( \eps_a \sinh s+ \frac{_{x \, \sinh s}}{^{\cosh s}}) \, (1+ \sin(\sqrt 2 (\eps_a s + x \, \tanh s))\, \cos (2 \theta)) 
\]
from which the formula for $\partial_x g(T_s, T_s) |_{x=0}$ follows at once.

Collecting these, we obtain the estimate for the mean curvature 
\[
H = {\mathcal O} \left( \tfrac{\log (1/\eps_a)  }{\cosh^2 s} \right) \, .
\]
This completes the estimate of the mean curvature. The estimate of its derivative follows similarly.  
\end{proof}

\section{The Analysis}

\subsection{Deformations of the Approximate Solution}
\mylabel{sec:deform}

The approximate solution $\tilde S_a$ constructed in the previous section is such that its mean curvature is almost equal to $H_a$ everywhere except in a small neighbourhood of each sub-lattice point, where it is nevertheless controlled by precise estimates.  The next task is to set up a means of finding a small deformation of $\tilde S_a$ whose mean curvature is exactly constant and equal to $H_a$.

To see how this can be done, let $f \in C^{2,\alpha}(\tilde S_a)$ be  given and let $N$ be a choice of unit normal vector field on $\papprox$.  Then if $f$ and its derivatives are sufficiently small (in a sense to be made precise in the next section), the neighbouring submanifold $\exp(f N) (\papprox)$ is an embedded submanifold of $S^3$ which is a small perturbation of $\papprox$.  Determining if $\exp(f N) (\papprox)$ has constant mean curvature is now a matter of solving a nonlinear partial differential equation.

\begin{defn} Define the \emph{deformation operator} to be the mapping $\Phi_a : C^{2,\alpha}(\papprox)  \rightarrow C^{0,\alpha}(\papprox)$ given by $\Phi_a (f) = H \big( \exp(f N)(\papprox) \big)$, where $H(\cdot)$ is the mean curvature operator.
\end{defn}

The deformation operator $\Phi_a$ is a non-linear partial differential operator acting on functions $f$ in $C^{2,\alpha}(\papprox)$ and so $\exp(f N) (\papprox)$ has constant mean curvature $H_a \in \R$ if and only if $f$ is a solution of the nonlinear partial differential equation $\Phi_a(f) = H_a$.  The linearization of the deformation operator at $0$ will also be needed.  The analysis of Section 2.1 asserts that $$\left. \frac{\dif}{\dif t} \right|_{t=0} \Phi_a ( tu) =  \Delta_a u + \Vert \tilde B_a \Vert^2 u + 2 u$$ where $\Delta_a$ is the Laplacian of $\papprox$ and $\tilde B_a$ is its second fundamental form.

\begin{defn}
    Denote $\mathcal{L}_a := \Dif \Phi_a (0)$.
\end{defn}

The tool that will be used to find the solution of the equation $\Phi_a (f) = H_a$ is the inverse function theorem.  See \cite{amr} for the proof of the version given here.

\begin{nonumthm}[IFT]
    Let $\Phi : \mathcal{B} \rightarrow \mathcal{B}'$ be a smooth map of Banach spaces, set $\Phi(0) = E$ and the denote the linearized operator $\Dif \Phi(0)$ by $\mathcal{L}$.  Suppose that $\mathcal{L}$ is bijective and the estimate $\Vert \mathcal{L} X \Vert \geq C \Vert X \Vert$ holds for all $X \in \mathcal{B}$.  Choose $R$ so that if $y \in \mathcal{B}$ is such that $\Vert Y \Vert \leq R$, then $\Vert \mathcal{L} X - \Dif \Phi(Y) X \Vert \leq \frac{1}{2} C \Vert X \Vert$.  If $ Z \in \mathcal{B}'$ is such that $\Vert Z - E \Vert \leq \frac{1}{2}C R$,  then there exists a unique $x \in \mathcal{B}$ with $\Vert X \Vert \leq R$ so that $\Phi(X) = Z$.  Moreover, $\Vert X \Vert \leq \frac{2}{C} \Vert Z - E \Vert$.
\end{nonumthm}

Finding the desired solution of the CMC equation by means of the inverse function theorem thus necessitates the following tasks.   First, appropriate Banach subspaces of $C^{2,\beta}(\papprox)$ and $C^{0,\beta}(\papprox)$ must be found --- along with appropriate norms --- so that the estimate of $\mathcal{L}_a$ can be achieved (this also establishes injectivity).  It must then be shown that $\mathcal{L}_a$ is surjective.  Next, estimates in these norms of the non-linear quantities (the size of $E := \Phi_a(0) - H_a$ and the size of the parameter $R$ giving the variation of $\Dif \Phi_a$) must be found.  Note that all these quantities depend \emph{a priori} on $a$.    Finally, the estimate of $E$ must be compared to the number $\frac{1}{2}C R$ and it must be shown that $\Vert E \Vert \leq \frac{1}{2} C R$ for all sufficiently small $a >0$.

\subsection{Function Spaces and Norms}

It is not possible to obtain a `good' linear estimate of the form $\Vert \mathcal{L}_a \, u \Vert \geq C\Vert u \Vert$ with any straightforward choice of Banach subspaces and norms.  There are essentially two reasons for this.  The first is that the operator $\mathcal{L}_a$ is \emph{not} injective on $C^{2,\alpha}(\papprox)$ due to the global Jacobi fields that come from $SO(4)$-rotations of the ambient $S^3$.  Each one-parameter family of rotations preserves the geometry of the ambient sphere --- and thus preserves the mean curvature of any submanifold of the sphere --- and so their generators are all elements of the kernel of $\mathcal{L}_a$.  The second reason for the absence of a good linear estimate is that $\mathcal{L}_a$ possesses \emph{small} eigenvalues so that even if one were to choose a Banach subspace transverse to the global Jacobi fields, the constant in the linear estimate would still depend on $a$ in an undesirable manner.  The eigenfunctions corresponding to these small eigenvalues come from \emph{approximate} Jacobi fields  such as the Jacobi fields that exist on the catenoidal necks of $\papprox$ and that `disappears' as $a \rightarrow 0$ and the necks pinch off.

The way in which the problems listed above will be dealt with here is twofold. First, the symmetries $\sigma_{j}$ and $\rho$ of the approximate solution must be exploited.  It turns out that the Jacobi fields, both approximate and true, do \emph{not} share these same symmetries.  Thus working in a space of functions possessing these symmetries will rule out the existence of small eigenvalues.  Second, it is necessary to use a somewhat non-standard norm to measure the `size' of functions in order to properly determine the dependence of the various estimates needed for the application of the inverse function theorem on the parameter $a$.  A \emph{weighted} H\"older norm will be used for this purpose, where each derivative term will be weighted by appropriate powers of a \emph{weight function} that accounts for the `natural' scaling property of the derivative operator.  The weight function is defined as follows.   Note that it is necessary only to precisely define the weight functions near a single sub-lattice point $\phi_0 (0,0)$ since the value of the weight function elsewhere can be found by symmetry.

\begin{defn}
    \mylabel{defn:weight}
    Define the \emph{weight function} $\zeta_a : \papprox \rightarrow \R$ in a neighbourhood of $\phi_0 (0,0)$ by
\begin{equation*}
    \zeta_a (x) =
    \begin{cases}
              \eps_a \, \cosh s  & \qquad \mbox{when} \:  x = \Xi (\eps_a \cosh s \cos \theta, \eps_a \cosh s \sin \theta , \eps_a \, s) \in \mathcal N_a \\ 
                \mbox{Interpolation} &\qquad \mbox{when} \: x \in {\mathcal T}_a \\
        |\mu|  &\qquad \mbox{when} \: x = \Xi (\mu ,t) \in {\mathcal C}^\pm_a 
          \end{cases}
\end{equation*}
\end{defn}

The required function spaces can now be defined.  First, recall the following notation.  If $\Omega \subseteq \papprox$ is any open subset of $\papprox$ and $q$ is any tensor on $\papprox$, then let
\begin{equation*}
    \Vert q \Vert_{0,\Omega} = \sup_{x \in \Omega} \Vert q(x) \Vert \\
    \qquad \mbox{and} \qquad [q]_{\alpha,\Omega} = \sup_{x,y \in \Omega} \frac{\Vert q(x) - \mathit{PT}(q(y))  \Vert}{\mathrm{dist}(x,y)^\alpha} \, ,
\end{equation*}
where the norms and the distance function that appear are taken with respect to the induced metric of $\papprox$, while $\mathit{PT}$ is the parallel transport operator from $x$ to $y$ with respect to this metric.   Now define the $C^{k,\alpha}_{\gamma}$ norm by
\begin{equation}
    \mylabel{eqn:weightnorm}
    \vert f \vert_{C^{k,\alpha}_{\gamma}(\Omega)} : = \vert \zeta_a^{-\gamma} f \vert_{0,\Omega} + \Vert \zeta_a^{-\gamma+1} \nabla f \Vert_{0,\Omega} + \cdots + \Vert \zeta_a^{-\gamma+k} \nabla^k f \Vert_{0,\Omega} + [ \zeta_a^{-\gamma+k+\alpha}\nabla^k f ]_{\alpha, \Omega} \, .
\end{equation}
Finally, let $\ckbg (\Omega)$ be the Banach space of $C^{k,\alpha}$ functions on $\Omega$ measured with respect to the norm \eqref{eqn:weightnorm}.

\begin{defn}
    The Banach spaces in which a solution of the deformation problem will be found are the spaces $\mathcal{B}^{k,\alpha,\gamma}_a : = \{ f \in C^{k,\alpha}_{\gamma} (\papprox) \: : \:   f \circ \sigma_{j} = f \circ \rho = f \: \: \forall \: j =1,2 \}$ of functions in $\ckbg(\papprox)$ possessing the symmetries $\sigma_{j}$ and $\rho$.  The parameters $\gamma$ will be chosen appropriately below.
\end{defn}

It remains to check that the operator $\Phi_a$ is well-behaved when acting on the Banach spaces $\Bkbg$.  It is straightforward to check that the map $\Phi_a : \ctbg(\papprox) \rightarrow \cobg (\papprox)$ is a smooth map in the Banach space sense and for any $u \in \Btbg$, the function $\Phi_a (u) : \papprox \rightarrow \R$ satisfies $\Phi_a (u) \circ \sigma_{j} = \Phi_a (u) \circ \rho = \Phi_a (u)$.  Thus $\Phi_a$ is well-defined as a map from $\Btbg$ to $\Bobg$.   The equivariance with respect to the symmetries $\sigma_{j}$ and $\rho$ is a consequence of the fact that these symmetries derive from isometries of the ambient Riemannian metric.  Finally, it is again straightforward to check that $\mathcal{L}_a$ is bounded in the operator norm by a constant independent of $a$, which is a consequence of the definition of the weighted norms being used here.

\subsection{The Linear Estimate}

The most important estimate in the solution of $\Phi_a (f) = H_a$ by means of the inverse function theorem is the estimate of the linearization $\mathcal{L}_a$ from below.  The purpose of this section is to prove this estimate.  The method used will be to construct an explicit solution of the equation $\mathcal{L}_a \, u = f$ by patching together local solutions on the neck region and away from the neck region.   This amounts to the construction of a right inverse for $\mathcal{L}_a$ --- which implies the surjectivity of $\mathcal{L}_a$, and by self-adjointness, the injectivity as well.

\begin{prop}
    \mylabel{thm:linest}
    Suppose $\gamma \in (-1,0)$. There exists $a_* >0$ such that, for all $a \in (0, a_*)$, the linearized operator $\mathcal{L}_a : \Btbg \rightarrow \Bobg$ satisfies the estimate
    $$ \vert u \vert_{\ctbg(\papprox)} \leq   C \, \eps_a^{\gamma} \,   \vert  \mathcal{L}_a \, u \vert_{\cobg(\papprox)}  $$
where $C$ is a constant independent of $a \in (0, a_*)$.
\end{prop}

\begin{proof}
    The patching argument requires a careful subdivision of $\papprox$ into separate pieces. First, we fix a parameter $\kappa >0$  smaller than $1/8$ of the least distance between two points of $\Lambda$ and define $\tilde {\mathcal C}^\pm_a (\kappa)$ to be the image of ${\mathbf R}^2 \setminus \bigcup_{\mu \in \Lambda} B_\kappa(\mu)$ by
    \[
    \mu \longmapsto \Xi (\mu, \pm (t_a - \eps_a \, \Gamma_\Lambda (\mu))
    \]
    and we define $\tilde {\mathcal N}_a(\kappa)$ to be the complement of $ \tilde {\mathcal C}^+_a (\kappa) \cup \tilde {\mathcal C}^-_a (\kappa) $ in $\tilde S_a$.
    
The next requirement is to carefully define two sets of cut-off functions with respect to these subdividing regions. All the following cut-off functions can be chosen smooth, as well as bounded by a constant independent of $a$ with respect to the $C^{k,\alpha}_0( \tilde S_a)$ norm.  Furthermore, each of these functions can be made symmetrical with respect to the symmetries $\sigma_j$ and $\rho$ satisfied by $\papprox$.  It is thus necessary to define them only in a neighbourhood of $\phi_0 (0,0)$ as follows.
\begin{itemize}

    \item $\chi_{\mathit{ext}, \kappa}^+$ (resp.~$\chi_{\mathit{ext}, \kappa}^-$) equals one in $\tilde {\mathcal C}^+ (2 \kappa)$ (resp.~$\tilde {\mathcal C}^- (2 \kappa)$) and equals zero in $\tilde {\mathcal N}_a (\kappa)\cup \tilde {\mathcal C}^-_a (\kappa)$ (resp. $\tilde {\mathcal N}_a (\kappa)\cup \tilde {\mathcal C}^+_a (\kappa)$).

    \item  $\chi_{neck, \kappa} : = 1 - \chi_{\mathit{ext},  \kappa}^+ - \chi_{\mathit{ext}, \kappa}^- $ and hence equals one n $\tilde {\mathcal N}_a (\kappa)$ and equals zero in $\tilde {\mathcal C}^\pm_a(2\kappa)$.
  
    \item $\eta_{\mathit{ext}}^+$ (resp.~$\eta_{\mathit{ext}}^-$) equals one in $\tilde {\mathcal C}^+ (2 \sqrt {\eps_a})$ (resp. $\tilde {\mathcal C}^- (2 \sqrt {\eps_a})$) and equals  zero in $\tilde {\mathcal N}_a (\sqrt {\eps_a})\cup \tilde {\mathcal C}^-_a (\sqrt {\eps_a})$ (resp.~$\tilde {\mathcal N}_a (\sqrt {\eps_a})\cup \tilde {\mathcal C}^+_a (\sqrt {\eps_a})$).

    \item  $\eta_{neck} : = 1 - \eta_{\mathit{ext}}^+ - \eta_{\mathit{ext}}^- $ and hence equals one in $\tilde {\mathcal N}_a (\sqrt {\eps_a})$ and equals zero in $\tilde {\mathcal C}^\pm_a(2\sqrt {\eps_a})$.

\end{itemize}
Observe that the cutoff functions  $\eta_{\mathit{ext}}^\pm$  and $\eta_{neck}$ can be chosen to be bounded in ${C}^{2, \alpha}_0 (\tilde S_a)$ uniformly in $a$.

To begin the patching argument, let $f \in \Btbg$ be given.  The idea is to construct an approximate solution of the equation ${\mathcal L}_a u = f$ by patching together a solution on the neck with a solution everywhere else.  This is carried out in the following two steps.

\paragraph*{Step 1.} Let $| \cdot |_{C^{k, \alpha}_\delta (\eps_a K)}$ denote the weighted $C^{k, \alpha}$ norm on $\eps_a \, K $, so that 
\[
| u |_{{C}^{k, \alpha}_\delta (\eps_a K)} : = |(\eps_a \, \cosh s)^{-\delta} u|_{0, \eps_a K} + \cdots + [ (\eps_a \, \cosh s)^{-\delta+ k + \alpha} \nabla^k u ]_{\alpha, \eps_a K}
\] 
where the norms and derivatives correspond to the metric on $\eps_a \, K$.  We are interested in functions that are invariant under the symmetry $u(s, \theta +\pi) = u(s, \theta)$.  The corresponding function spaces will be denoted by ${C}^{k, \alpha}_{\delta , sym} (\eps K)$. We have the following result (whose proof follows from a similar proof that can be found for example in \cite{kusnermazzeopollack} in the context of constant mean curvature surfaces) and which follows from the result of Lemma~\ref{lemma:catenoidjac}. The operator 
\[
\begin{array}{rcccllll}
L_{\eps_a K} :  & {C}^{2, \alpha}_{\delta , sym} (\eps_a K)   & \longrightarrow & {C}^{0, \alpha}_{\delta -2, sym} (\eps_a K) \\[3mm]
& u & \longmapsto & {\mathcal L}_{\eps_a K}\, u
\end{array}
\]
is injective if $\delta < 0$ and hence it is surjective if $\delta  > 0$, $\delta \notin {\mathbf N}$.  This latter fact follows from a duality argument in weighted Lebesgue spaces and elliptic regularity theory can be used to prove the corresponding result in weighed H\"older spaces.

We define
\[
\tilde f_0^{(1)} (s, \theta) : = Ê(2 \chi (s) -1) \, \tanh s \qquad \mbox{and} \qquad \tilde f_0^{(2)} (s, \theta) : = Ê(2 \chi (s) -1) \, (s \,  \tanh s -1)
\] 
where $\chi$ is a cut-off function identically equal to $0$ for $s < -1$ and identically equal to $1$ for $s > 1$ and we recall that we have already defined  in Lemma~\ref{lemma:catenoidjac}
\[
f_0^{(1)} (s, \theta) : = Ê\tanh s \qquad \mbox{and} \qquad f_0^{(2)} (s, \theta) : = Ê s \,  \tanh s -1 \, .
\] 
The deficiency space  ${\mathcal D}$ is defined to be the space of functions 
\[
{\mathcal D} : =  \mbox{Span} \, \{Êf_0^{(1)} , \tilde f_0^{(1)} ,  f_0^{(2)} , \tilde f_0^{(2)}     \}  \, .
\]
Then, for all $\delta \in (-2, 0)$,  the operator
\[
\begin{array}{rcccllll}
\tilde L_{\eps_a K} :  & {C}^{2, \alpha}_{\delta , sym} (\eps_a K) \oplus {\mathcal D} & \longrightarrow & {C}^{0, \alpha}_{\delta -2, sym} (\eps_a K) \\[3mm]
& u & \longmapsto & {\mathcal L}_{\eps_a K}\, u
\end{array}
\]
is surjective and has a two dimensional kernel  \cite{kusnermazzeopollack}.  This result follows from the fact that if $f \in {C}^{0, \alpha}_{\delta -2, sym} (\eps_a K)$ and $\delta \in (-2,0)$, then $f \in {C}^{0, \alpha}_{-\delta -2, sym} (\eps_a K)$ and hence one can find $u \in {C}^{2, \alpha}_{-\delta , sym} (\eps_a K)$ solution of  ${\mathcal L}_{\eps_a K} u =f$. Then one checks that, for  $s >0$, the function $u$ can be decomposed into
\[
u(s, \theta ) =  v(s, \theta) + a \, f_0^{(1)} + b \,f_0^{(2)} 
\]
where $v$ is bounded by a constant times $(\cosh s)^{\delta}$ for $s >0$ and $a, b \in {\mathbf R}$. A similar decomposition is available for $s <0$. 

Observe that the kernel of $\tilde L_{\eps_a K}$ is simply generated by the functions $f_0^{(1)}$ and $f_0^{(2)} $ that are defined in Lemma~\ref{lemma:catenoidjac}. Moreover, one can find a right inverse whose norm does not depend on $\eps_a$. For example it is possible to choose a right inverse that maps into $ {C}^{2, \alpha}_{\delta , sym} (\eps_a K) \oplus {\mathcal D}_0$ where ${\mathcal D}_0 : =  \mbox{Span} \, \{Ê\tilde f_0^{(1)}  , \tilde f_0^{(2)}  \}$.

We fix $\kappa_0$ small enough and define 
\[
f_{neck} : =  \chi_{neck, \kappa_0} \, f
\]
and apply the above result when $\delta = \gamma \in (-2,0)$.  Observe that $|f_{neck}|_{\cobg (\eps K) }Ê\leq C \, | f |_{\cobg (\tilde S_a) }$. Therefore, we can define $u_{neck} \in {C}^{2, \alpha}_{\gamma , sym} (\eps_a K) \oplus {\mathcal D}_0$  solution of ${\mathcal L}_{\eps_a K } \, u_{neck} = f_{neck}$. It is possible to decompose 
\[
u_{neck} = v_{neck} + p_1  \, \tilde f_0^{(1)}   + p_2 \,  \tilde f_0^{(2)} 
\]
where $v_{neck} \in \ctbg (\eps_a K )$ and $p_j$ are constants.  Furthermore, one has the estimate
\begin{equation}
    \label{eqn:neckdecompest}
    | v_{neck} |_{\ctbg (\eps_a K) }  + \eps_a^{-\gamma} (|p_1| + |p_2| )\leq C  \, | f_{neck}|_{\cobg (\eps_a K)} \, .
   \end{equation}
We extend this solution to all $\tilde S_a$ as follows
\[
\bar u_{neck} : = \chi_{neck} \,  v_{neck} + \eta_{neck} (p_1  \, \tilde f_0^{(1)}   + p_2 \,  \tilde f_0^{(2)} + q_1 \, f_0^{(1)} + q_2 \, f_0^{(2)}) + \eta_{\mathit{ext}}^+ \,  r^+ \, \Gamma_\Lambda  +\eta_{\mathit{ext}}^- \,  r^- \, \Gamma_\Lambda
\]
where the coefficients $q_j$ and $r_j$ are  determined to ensure a good matching of the functions $r^\pm  \, \Gamma_\Lambda$ and $p_1  \, \tilde f_0^{(1)}   + p_2 \,  \tilde f_0^{(2)} + q_1 \, f_0^{(1)} + q_2 \, f_0^{(2)} $ on the different summands. Namely, we find that these coefficients must satisfy the system  
\begin{equation}
	\begin{gathered} 
		r ^+ \, (\gamma_\Lambda - \log (\eps_a/2)  ) - q_1 + q_2 =  p_1 - p_2 \\
		r ^- \, (\gamma_\Lambda - \log (\eps_a/2)  ) + q_1 + q_2 =  p_1 + p_2 \\
		r^+  + q_2  =  -  p_2  \\
		r^-  + q_2 = p_2 \, .
	\end{gathered}
\label{eq:s1}
\end{equation}
For example, when $s \sim \frac{1}{2}Ê\, \log 1/ \eps_a$, we can write  
\[
p_1  \, \tilde f_0^{(1)}   + p_2 \,  \tilde f_0^{(2)} + q_1 \, f_0^{(1)} + q_2 \, f_0^{(2)}  \sim (p_1+q_1) + (p_2+q_2) \, (s -1)
\]
while, when $r \sim \sqrt \eps_a$ we can write 
\[
r^+  \, \Gamma_\Lambda \sim r^+  \, ( - \log (\eps_a/2) - s + \gamma_\Lambda )
\]
if we change variable $|\mu| = \eps_a \cosh s$ with $s \sim  \frac{1}{2}Ê\, \log 1/ \eps_a$. The first and third identities in (\ref{eq:s1}) are obtained by equating that coefficients of the constant function and the function $s$ in these two expansions. 

Therefore, we have the estimates 
\[
(\log 1/\eps_a)^{-1} \, | q_1| +(\log 1/\eps_a) \, | q_2| +  | r^\pm| \leq C \, (|p_1|+ |p_2|) \leq C \, \eps_a^\gamma \, | f_{neck}|_{{ C}^{0, \alpha}_{\gamma-2} (\eps_a K)}.
\]
Putting all of this together, we find the estimate
\[
| \bar u_{neck}  |_{{C}^{2, \alpha}_{\gamma} (\tilde S_a)} \leq C \,  \eps_a^{\gamma}  \, |  f |_{{C}^{0, \alpha}_{\gamma-2} (\tilde S_a)} \, .
\]

Now using the fact that the Jacobi operator in $\tilde {\mathcal N}_a$ is close to the Jacobi operator on $\eps_a \, K$, we can evaluate  ${\mathcal L}_{a} \, \bar u_{neck} - f_{neck}$. For all  $\kappa_1 \leq  \kappa_0 / 4$, we find
\begin{equation}
|  \chi_{neck ,\kappa_1} \, ({\mathcal L}_{a} \, \bar u_{neck} - f ) |_{{ C}^{0, \alpha}_{\gamma-2} (\tilde S_a)} \leq C \, ( \kappa_1^2 + \eps_a^{1+\gamma } \ (\log 1/\eps_a)  ) \, |  f |_{{ C}^{0, \alpha}_{\gamma-2} (\tilde S_a)}
\label{eq:estimate1}
\end{equation}
where the constant $C >0$ does not depend on $\kappa_1$. To obtain this estimate, it is enough to estimate the difference between the Jacobi operators ${\mathcal L}_{\eps_a K}$ and ${\mathcal L}_a$. For example, using the analysis of the proof of Proposition~\ref{prop:estH} we find that the metric of ${\mathcal N}_a$ is not far from the metric of $\eps_a K$ and indeed 
\[
g_{a} = \frac{_1}{^2} \, \eps_a^2 \, \cosh^2 s \, (ds^2 + d\theta^2) + {\mathcal O} (\eps_a^3 \, (\log 1/\eps_a) \, \cosh^2 s )
\]
if $\Xi (\eps_a \cosh s \cos \theta, \eps_a \cosh s \cos \theta , \eps_a s ) \in {\mathcal N}_a (\kappa)$ while 
\[
g_{\eps_a K} = \frac{_1}{^2} \, \eps_a^2 \, \cosh^2 s \, (ds^2 + d\theta^2)
\]
Similarly, the square of the norm of the second fundamental form of $\tilde {\mathcal N}_a$  is again not far from the square of the norm of the second fundamental form on  on $\eps_a K$ and indeed 
\[
|B_{a}|^2 = \frac{_4}{^{\eps_a^2 \, \cosh^4 s}} +  {\mathcal O} (\frac{_{\log 1/\eps_a}}{^{\eps_a \, \cosh^4 s}} ) + {\mathcal O} (\frac{_1}{^{\eps_a \, \cosh^2 s}} )
\]
in ${\mathcal N}_a (\sqrt{\eps_a})$
and 
\[
|B_{a}|^2 =  2 + {\mathcal O} (\eps_a \, (\log 1/\eps_a))
\]
in ${\mathcal C}_a^\pm (\sqrt{\eps_a})$ while 
\[
| B_{\eps_a K} |^2 = \frac{_4}{^{\eps_a^2 \, \cosh^4 s}} 
\]
It is then enough to use the fact that ${\mathcal L}_{\eps_a K} = \Delta_{g_{\eps_a K}} + |ÊB_{\eps_a K}|^2$ while ${\mathcal L}_{a} = \Delta_{a} + |ÊB_{a}|^2 +2$.  The contribution to (\ref{eq:estimate1}) of the difference between ${\mathcal L}_a \bar u_{neck}$ and $f$  in ${\mathcal N}_a (\sqrt{\eps_a})$ can be estimated by a constant times $\eps_a^{1+\gamma } \ (\log 1/\eps_a)  ) \, |  f |_{{ C}^{0, \alpha}_{\gamma-2} (\tilde S_a)}$.   Using the fact that $\Gamma_\Lambda$ is annihilated by ${\mathcal L}_{C_0}$, we find that the contribution to (\ref{eq:estimate1}) of the difference between ${\mathcal L}_a \bar u_{neck}$ and $f$  in ${\mathcal C}_a^\pm (\sqrt{\eps_a})\setminus {\mathcal C}_a^\pm (\kappa_1)$ can be estimated by a constant times $( \kappa_1^2 + \eps_a^{1+\gamma}  \ (\log 1/\eps_a)  ) \, |  f |_{{ C}^{0, \alpha}_{\gamma-2} (\tilde S_a)}$. Finally, the influence of the cutoff functions $\eta_{neck}$ and $\eta_{ext}^\pm$ in ${\mathcal T}_a$ produces a discrepancy that can be evaluated by a constant times $\eps_a^{1+\frac{\gamma}{2} } \ (\log 1/\eps_a)  \, |  f |_{{ C}^{0, \alpha}_{\gamma-2} (\tilde S_a)}$ and hence of a much smaller magnitude.

Observe that, when $\gamma \in (-1, 0)$, we have 
\[
|  {\mathcal L}_{a} \, \bar u_{neck} - f |_{{ C}^{0, \alpha}_{\gamma-2} (\tilde S_a)} \leq C \, |  f |_{{ C}^{0, \alpha}_{\gamma-2} (\tilde S_a)}
\]
provided $a$ is small enough. This follows at once from (\ref{eq:estimate1}) together with the fact that $\Gamma_\Lambda$ is annihilated by ${\mathcal L}_{C_0}$ away from the points of $\phi_0(\Lambda)$.

Finally, we set
\[
\hat  f_{\mathit{ext}}^\pm  : =  \chi_{ext, \kappa_1}^\pm  \, \left( f -  {\mathcal L}_{a} \,  \bar u_{neck}  \right) .
\]
The functions $\hat  f_{\mathit{ext}}^\pm$ being supported away from the neck, they can be considered as functions on $\tilde S_a$ or functions on ${\mathcal C}^\pm_a$ or even functions on $C_0$. Observe that 
\[
|  \hat  f_{\mathit{ext}}^\pm   |_{{ C}^{0, \alpha} (C_0)} \leq C_{\kappa_1}  \, |  f |_{{ C}^{0, \alpha}_{\gamma-2} (\tilde S_a)}
\]
for some constant $C_{\kappa_1}$ that depends on $\kappa_1$ and $\gamma$.

\paragraph*{Step 2.} Let $| \cdot |_{{C}^{k, \alpha}_\delta ( C_0  \setminus \phi_0(\Lambda ) ) }$ denote the weighted ${C}^{k, \alpha}$ norm on $ C_0  \setminus \phi_0(\Lambda ) $, so that 
\[
| \cdot |_{{C}^{k, \alpha}_\delta ( C_0  \setminus \phi_0 (\Lambda ) ) } : = ||\mu|^{-\delta} u|_{0, \Omega} + \cdots + [ |\mu|^{-\delta+ k + \alpha} \nabla^k u ]_{\alpha, \Omega}
\]
where $\Omega$ is a fundamental cell of the lattice $\phi_0(\Lambda )$ containing $\phi_0(0,0)$.

We first find $u_{\mathit{ext}}^\pm$ the unique $C^{2,\alpha}(C_0)$ solution of  ${\mathcal L}_{C_0} u_{\mathit{ext}}^\pm = \hat f_{\mathit{ext}}^\pm $ on $ C_0$ satisfying the usual invariance property (\ref{eq:eqi}).  Observe that, near $\mu =(0,0)$ the Taylor expansion of this solution is given by
\begin{equation}
    \label{eqn:extexp}
    u_{\mathit{ext}}^\pm = u_{\mathit{ext}}^\pm (0)  + v_{\mathit{ext}}^\pm   
\end{equation}
where  $v_{\mathit{ext}}^\pm \in C^{2, \alpha}_{2}(C_{0}\setminus \phi_0 (\Lambda) )$.  This reflects the fact that $u_{\mathit{ext}}^\pm $ satisfies $(\Delta + 2) u_{\mathit{ext}}^\pm (0)  =0$ near $(0,0)$ and hence $u_{\mathit{ext}}^\pm $ is smooth there. Moreover, since $u_{\mathit{ext}}^\pm (-\mu)= u_{\mathit{ext}}^\pm (\mu) $, then the first partial derivatives of $u_{\mathit{ext}}^\pm$ vanish at $(0,0)$. We also have the estimate
\begin{equation}
    \label{eqn:extexpest}
    | v_{\mathit{ext}}^\pm |_{{C}^{2, \alpha}_{2} (C_0 \setminus \phi_0(\Lambda)) } + | u_{\mathit{ext}}^\pm (0)  |  \leq C \, | \hat f_{\mathit{ext}}^\pm |_{{C}^{0, \alpha} ( C_0)} \leq C_{\kappa_1} \, | f  |_{{C}^{0, \alpha}_{\gamma-2} (\tilde S_a)}
\end{equation}
holds for some constant $C$ independent of $f$ and $a$ (but $C_{\kappa_1} $ depends on $\kappa_1$).  As in Step 1, we extend $u_{\mathit{ext}}^\pm $ to $\tilde S_a$ using the function $\Gamma_\Lambda$ and the functions $f_0^{(j)}$. We define 
\begin{align*}
	\bar u_{ext} := \eta_{\mathit{ext}}^+ \,  (     u_{\mathit{ext}}^+  +      \hat r^+ \, \Gamma_\Lambda )  +\eta_{\mathit{ext}}^- \, (u_{\mathit{ext}}^-   +  \hat r^- \, \Gamma_\Lambda )	 + \eta_{neck} (\hat q_1 \, f_0^{(1)} + \hat q_2 \, f_0^{(2)})
	 \end{align*}
where the coefficients $\hat q_j$ and $\hat r_j$ are  determined to ensure a good matching of the functions $\hat q_1 \, f_0^{(1)} + \hat q_2 \, f_0^{(2)}$ and $  u_{\mathit{ext}}^\pm (0) + \hat r^\pm  \, \Gamma_\Lambda$ on the different summands. Namely, we find that these coefficients must satisfy the system  
\begin{equation*}
	\begin{gathered}
		\hat r ^+ \, (\gamma_\Lambda - \log (\eps_a/2)  )  + u^+_{ext} (0) =  \hat q_1 - \hat q_2 \\
		\hat r ^- \, (\gamma_\Lambda - \log (\eps_a/2)  )  + u_{ext}^- (0) =  - \hat q_1  - \hat q_2 \\
		\hat r^+ =  -  \hat q_2  \\
		\hat r^-  =  - \hat q_2 .
	\end{gathered}
\end{equation*}
Therefore, we have the estimates 
\[
(\log 1/\eps_a) \, ( | \hat q_2| + | r^\pm| ) + |\hat q_1|Ê\leq C \, ( | \hat f_{ext}^+ |_{{C}^{0, \alpha}} +| \hat f_{ext}^- |_{{C}^{0, \alpha}} ) \leq C_{\kappa_1} \, | f  |_{{C}^{0, \alpha}_{\gamma-2} (\tilde S_a)} \, .
\] 
Putting all of this together, we obtain the estimate
\[
|  \bar u_{ext}  |_{{C}^{2, \alpha}_{\gamma} (\tilde S_a)} \leq C_{\kappa_1} \, |  f |_{{C}^{0, \alpha}_{\gamma-2} (\tilde S_a)} \, .
\]

Now using the fact that the Jacobi operator on $ C_0$ is close to the Jacobi operator on $\tilde {\mathcal C}_a^\pm$, we can evaluate  ${\mathcal L}_{a} \, \bar u_{ext} - \hat f_{ext}^+ -\hat f_{ext}^- $.  With little work, and using the strategy developed in step 1, we find
\begin{equation}
|  {\mathcal L}_{a} \, \bar u_{ext} - \hat f_{ext}^+ - \hat f_{ext}^-  |_{{C}^{0, \alpha}_{\gamma-2} (\tilde S_a)} \leq C_{\kappa_1}  \, \eps_a \, (\log 1/\eps_a) \,  |  f |_{{C}^{0, \alpha}_{\gamma-2} (\tilde S_a)} 
\label{eq:estimate2}
\end{equation}
There is no difficulty in obtaining this estimate. Observe that we have used the fact that $\Gamma_\Lambda$ is annihilated by ${\mathcal L}_{C_0}$ and also the fact that the functions $f^{(j)}_0$ are in the kernel of ${\mathcal L}_{\eps_a K}$.
The contribution to (\ref{eq:estimate2}) of the difference between ${\mathcal L}_a \bar u_{ext}$ and $ \hat f_{ext}^+ + \hat f_{ext}^-$  in ${\mathcal N}_a (\sqrt{\eps_a})$ can be estimated by a constant (depending on $\kappa_1$) times $\eps_a \ (\log 1/\eps_a)  \, |  f |_{{ C}^{0, \alpha}_{\gamma-2} (\tilde S_a)}$.   The contribution to (\ref{eq:estimate2}) of the difference between ${\mathcal L}_a \bar u_{ext}$ and $\hat f_{ext}^+ + \hat f_{ext}^-$  in ${\mathcal C}_a^\pm (\sqrt{\eps_a})\setminus {\mathcal C}_a^\pm (\kappa_1)$ can be estimated by a constant  (depending on $\kappa_1$) times $ \eps_a^{2 - \frac{\gamma}{2}}  \ (\log 1/\eps_a)   \, |  f |_{{ C}^{0, \alpha}_{\gamma-2} (\tilde S_a)}$. Finally, the influence of the cutoff functions $\eta_{neck}$ and $\eta_{ext}^\pm$ in ${\mathcal T}_a$ produces a discrepancy that can be evaluated by a constant  (depending on $\kappa_1$) times $\eps_a^{1-\frac{\gamma}{2} }  \, |  f |_{{ C}^{0, \alpha}_{\gamma-2} (\tilde S_a)}$ and hence of a much smaller magnitude.

Collecting the estimates we conclude that 
\[
|  \bar u_{ext} + \bar u_{neck} |_{{ C}^{2, \alpha}_{\gamma} (\tilde S_a)} \leq C_{\kappa_1}  \,  \eps_a^{\gamma}  \, |  f |_{{ C}^{0, \alpha}_{\gamma-2} (\tilde S_a)}
\]
and also that 
\[
|  {\mathcal L}_{a} \, (\bar u_{ext} + \bar u_{neck}) - f  |_{{ C}^{0, \alpha}_{\gamma-2} (\tilde S_a)} \leq  ( C\, \kappa_1^2 + C_{\kappa_1} \,  \eps_a^{1+\gamma} \, (\log 1/\eps_a) )  \, |  f |_{{ C}^{0, \alpha}_{\gamma-2} (\tilde S_a)}
\]
In other words, for all $a$ small enough, the mapping $ f \longmapsto \bar u_{ext} + \bar u_{neck}$ is almost a right inverse for ${\mathcal L}_a$ and the result now follows from a standard perturbation argument, provided $\kappa_1$ is fixed small enough and $\gamma \in (-1, 0)$.
\end{proof}

\subsection{The Nonlinear Estimates and the Conclusion of the Proof}

\paragraph*{Nonlinear Estimates.}  As mentioned earlier, the proof of the Main Theorem requires two more estimates.  First, it is necessary to show that $\left\vert \Phi_a (0) - H_a \right\vert$ is small in the $\cobg (\papprox)$ norm.  Second, it is necessary to show that $\left[ \Dif \Phi_a (f)  - \mathcal{L}_a \right]$ can be made to have small operator norm (with respect to the $\ckbg (\papprox)$ norms) if $f$ is chosen sufficiently small in the $\ctbg (\papprox)$ norm.  Note that $\gamma \in (-1,0)$.  Once these estimates are understood, it will be possible to conclude the proof of the Main Theorem simply by invoking the inverse function theorem

The following is  a simple consequence of the result of Proposition \ref{prop:estH}.
\begin{prop}
    \mylabel{prop:error}
    The quantity $\Phi_a (0)$, which is the mean curvature of $\papprox$, satisfies the following estimate.  If $a$ is sufficiently small, then there exists a constant $C_*$ independent of $a$ so that
    \begin{equation}
        \mylabel{eqn:error}
        \big\vert \Phi_a (0) - H_a  \big\vert_{\cobg} \leq C_* \, \eps_a^2 \, (\log 1/\eps_a).
    \end{equation}
\end{prop}

\noindent We also need the following.
\begin{prop}
    \mylabel{prop:nonlin}
    The differential of the operator $\Phi_a$ satisfies the following estimate.  Fix $\kappa_2  >0$. Then the exists a constant $C_{\kappa_2} >0$ such that, for all  sufficiently small $a$ and all $f$ satisfying $\vert f \vert_{\ctbg} \leq  \kappa_2 \, \eps_a^{2+\gamma} \, \log 1/\eps_a$ we have
    \begin{equation}
        \mylabel{eqn:nonlin}
        \big\vert\Dif \Phi_a (f) \, u - \mathcal{L}_a  \, u \big\vert_{\cobg (\papprox)} \leq C_{\kappa_2} \, \eps_a^{1+2 \gamma} \, (\log 1/\eps_a)Ê\,  \vert u \vert_{\ctbg (\papprox)}
    \end{equation}
     for any $u \in \ctbg (\papprox)$.
\end{prop}

\begin{proof}
The derivative of $\Phi_a$ at a point different from $0$ is a similar calculation as at the point $0$.  That is
$$\left. \frac{\dif \Phi_a (f + t u)}{\dif t} \right|_{t=0} = \Delta_{a, f} u + \big(\Vert \tilde B_{a, f}\Vert^2 + 2 \big) u$$
where $\Delta_{a, f}$ and $\tilde B_{a, f}$ are the Laplacian and the second fundamental form, respectively,  of the submanifold $\tilde S_{a, f} : =  \exp(fN)(\papprox)$.  Consequently, we must estimate the quantity
\begin{equation}
    \mylabel{eqn:nonlinsetup}
 \big( \Dif \Phi_a (f) \, u - \mathcal{L}_a \, u  \big)  =  \big( \Delta_{a,f} - \Delta_a \big) u +    \big( \Vert \tilde B_{a ,f} \Vert^2- \Vert \tilde B_a \Vert^2 \big) 
\end{equation}
in the  $\vert \cdot  \vert_{C^{0, \alpha}_{\gamma -2} (\papprox)}$ norm.   To obtain this estimate, it is sufficient to estimate the difference between the metrics on $\tilde S_{a,f}$ and $\tilde S_a$ which will provide the estimates between the Laplacians and also the difference between the square norms of the second fundamental forms of these surfaces. 

Observe that on the normal graph of $f$ over ${\mathcal C}^\pm_a(2{\eps_a})$, the surface $\tilde S_{a,f}$ can be represented as graphs over $C_0 \setminus \bigcup_{p\in \phi_0(\Lambda)} B_{2 \sqrt{\eps_a}} (p)$ for some functions $\pm u_a + \tilde f^\pm$ where $\tilde f^\pm$ satisfies
\[
| \chi_{ext, 2\eps_a}^\pm \, \tilde f^\pm |_{C^{2, \alpha}_\gamma (C_0 \setminus \phi_0(\Lambda))} \leq C_{ \kappa_2} \, \eps_a^{2+\gamma} \, \log 1/\eps_a
\]
where the cutoff functions $\chi_{ext, 2\eps_a}^\pm$ are assumed to be chosen so that their norm in $C^{2, \alpha}_0 (\tilde S_a)$ are bounded independently of $a$. It follows at once from the analysis performed in the proof of Proposition~\ref{prop:estH} that
\[
| \chi_{ext, 2\eps_a}^\pm \, (\Delta_{a,f} - \Delta_a ) u  |_{C^{0, \alpha}_{\gamma-2} (C_0 \setminus \phi_0(\Lambda))} \leq C_{ \kappa_2} \, \eps_a^{1+2\gamma} \, (\log 1/\eps_a) \, | \chi_{ext, 2\eps_a}^\pm \, u |_{C^{2, \alpha}_{\gamma} (C_0 \setminus \phi_0(\Lambda))} 
\]
and also that 
\[
| \chi_{ext, 2\eps_a}^\pm \, \big( \Vert \tilde B_{a ,f} \Vert^2- \Vert \tilde  B_a \Vert^2 \big)  \, u  |_{C^{0, \alpha}_{\gamma-2} (C_0 \setminus \phi_0(\Lambda))} \leq C_{ \kappa_2} \, \eps_a^{1+2\gamma} \, (\log 1/\eps_a) \, | \chi_{ext, 2\eps_a}^\pm \, u  |_{C^{2, \alpha}_{\gamma} (C_0 \setminus \phi_0(\Lambda))} 
\]
Let us explain where these estimates come from. The first estimate follows from the fact that the coefficients of the metric of the graph of a function $u$ contains terms of the form  $\partial_i u \,\partial_j u$ which produce a discrepancy between the coefficients of the metric of $\tilde S_{a,f}$ and $\tilde S_a$ that involves terms of the form $\partial_i u_{a} \, \partial_j \tilde f^\pm$.   The second estimate follows from the fact that the difference between the coefficients of the square of the fundamental form of the graph of a function $u$ contains  terms of the form  $\partial_{ij} u \, \partial_{kl} u$ which produce a discrepancy between the square of the fundamental forms of  $\tilde S_{a,f}$ and $\tilde S_a$ that involves terms of the form $\partial_{ij} u_{a} \, \partial_{kl} \tilde f$. 

The estimate in ${\mathcal N}_a (2 \eps_a)$ follows easily from the analysis of the second half of the proof of Proposition~\ref{prop:estH} and does not yield worse estimates than we already have. We leave the details to the reader.
\end{proof}

\paragraph*{Conclusion of the Proof.} The estimates for the proof of the Main Theorem are now all in place and the conclusion of the theorem becomes a simple verification of the conditions of the inverse function theorem. We choose $\gamma \in (-1/3, 0)$.  First by Theorem \ref{thm:linest}, the linearization satisfies the estimate
\[
\vert u \vert_{\ctbg(\papprox)}  \geq C^* \, \eps_a^{-\gamma}  \, \vert \mathcal{L}_a \, u \vert_{\cobg(\papprox)}  
\]
where $C^*$ is a constant independent of $a$.  Therefore by the inverse function theorem of Section \ref{sec:deform} along with Proposition \ref{prop:nonlin}, a solution of the deformation problem can be found if $$\big\vert \Phi_a (0) - H_a \big\vert_{\cobg(\papprox)} \leq \frac{1}{2} C^* \, \eps_a^{-\gamma}  \, R$$
where $R = \kappa_2 \, \eps_a^{2+\gamma} \, (\log 1/\eps_a) $ and if 
\[
 \big\vert\Dif \Phi_a (f) \, u - \mathcal{L}_a  \, u \big\vert_{\cobg (\papprox)} \leq \frac{1}{2} C^* \, \eps_a^{-\gamma}  \,  \vert u \vert_{\ctbg (\papprox)} \, .
\]
But the above result shows that $ \big\vert\Dif \Phi_a (f) \, u - \mathcal{L}_a  \, u \big\vert_{\cobg (\papprox)} \leq C_{\kappa_2} \, \eps_a^{1+2 \gamma} \, (\log 1/\eps_a)Ê\,  \vert u \vert_{\ctbg (\papprox)}$ and Proposition \ref{prop:error} shows that $\big\vert\Phi_a (0) -  H_a \big\vert_{\cobg(\papprox)} \leq   C_* \, \eps_a^{2} \, (\log 1/\eps_a)$.  Hence this can always be done if $a$ is sufficiently small and $\kappa_2$ is large enough to ensure $\kappa_2 \, C^*  \geq C_*$.  This concludes the proof of the Main Theorem. \hfill \qedsymbol

\vfill
\begin{multicols}{2}
\parindent = 0ex
\renewcommand{\baselinestretch}{1.0}
\normalsize

Adrian Butscher

University of Toronto

1265 Military Trail

Toronto, Ontario, M1C 1A4

Canada

Frank Pacard

Universit\'e de Paris XII

61 rue du G\'en\'eral Charles de Gaulle

94010 Cr\'eteil

France

\end{multicols}

\newpage
\renewcommand{\baselinestretch}{1.0}
\normalsize
\bigskip \bigskip

\bibliography{cmc}
\bibliographystyle{amsplain}

\end{document}